\def\DocTitle{Sparse approximations for contact mechanics}
\title{\DocTitle}
\ifx\fmt\undefined
    \documentclass{article}
    \usepackage[a4paper]{geometry}
    \newenvironment{keywords}{\textbf{Keywords:}}{}
\else
    \documentclass[AMA,Times1COL]{WileyNJDv5}
    \articletype{Research Article}%
\fi

\usepackage{graphicx}
\usepackage{amsfonts}
\usepackage{amsmath,amsthm}
\usepackage{fleqn}
\usepackage{epstopdf}
\usepackage{float}
\usepackage{array}
\usepackage{xcolor,colortbl}

\usepackage{booktabs}
\usepackage{multirow}
\usepackage{arydshln}
\usepackage{enumerate}

\ifx\fmt\undefined
    \newtheorem{remark}{Remark}
    \theoremstyle{definition}
    \newtheorem{hypothesis}{Hypothesis}
\fi

\definecolor{ncred}{rgb}{0.88,0,0}
\usepackage{hyperref}
\hypersetup{
    pdftitle={Sparse approximations for contact mechanics},
    pdfauthor={Kiran Sagar Kollepara},
    pdfsubject={Research article},
    pdfkeywords={contact mechanics, variational inequalities, reduced order models, sparse decomposition, dictionary based approximation},
    bookmarks=true,
    bookmarksopen=true,
    bookmarksopenlevel=1,
    bookmarksnumbered=true,
    pdfpagemode=UseOutlines,
    colorlinks=true,
    linkcolor=blue,
    citecolor=ncred,
}

\makeatletter
\newcounter{HALG@line}
\renewcommand{\theHALG@line}{\thealgorithm.\arabic{ALG@line}}
\makeatother

\usepackage[capitalize]{cleveref}
\crefname{hypothesis}{Hypothesis}{Hypotheses}

\crefformat{algorithm}{Alg.~#2#1#3}
\crefformat{line}{Step~#2#1#3}
\crefrangeformat{line}{Steps~#3#1#4 to~#5#2#6}
\crefmultiformat{line}{Steps~#2#1#3}{ and~#2#1#3}{, #2#1#3}{ and~#2#1#3}

\usepackage{subcaption}
\usepackage{bm}
\usepackage{algpseudocode}
\usepackage{algorithm}
\usepackage{todonotes}
\usepackage{mathtools}

\ifx\fmt\undefined
    \usepackage[square,numbers,sort&compress]{natbib}
    \usepackage{doi}
\fi

\usepackage[final]{showlabels} %% to be removed (or include [final] option) in the end 

\DeclareMathOperator*{\argmin}{arg\,min}
\DeclareMathOperator*{\argmax}{arg\,max}

\newcommand{\norm}[2]{{\left| \left| #1 \right| \right|}_{#2}}

\newcommand{\PrimalRB}[0]{\mathbf{\Phi}}

\newcommand{\DualDict}[0]{\mathbf{D}_d}
\newcommand{\DualDictActive}[0]{\mathbf{D}_{d,\mathcal{I}}}
\newcommand{\PrimalDict}[0]{\mathbf{D}_p}

\newcommand{\vh}[0]{{\bm{v}^h}}
\newcommand{\uh}[0]{{\bm{u}^h}}

\newcommand{\lambdah}[0]{{\bm{\lambda}^h}}

\begin{document}
\ifx\fmt\undefined
    \author{Kiran Sagar Kollepara
        \and
        Jos\'e V. Aguado
        \and
        Yves Le Guennec
        \and
        Luisa Silva
        \and
        Domenico Borzacchiello
            }
\else
    \author[1]{Kiran Sagar Kollepara}
    \author[1]{Jos\'e V. Aguado}
    \author[2]{Yves Le Guennec}
    \author[1]{Luisa Silva}
    \author[1]{Domenico Borzacchiello}

    \authormark{KOLLEPARA \textsc{et al.}}
    \titlemark{\DocTitle}

    \address[1]{\orgdiv{Institut de Recherche en Génie Civil et Mécanique (GeM)}, \orgname{{E}cole Centrale de Nantes}, \orgaddress{\state{Nantes}, \country{France}}}
    \address[2]{\orgname{IRT Jules Verne}, \orgaddress{\state{Bouguenais}, \country{France}}}

    \corres{Corresponding author Kiran Sagar Kollepara\\ \email{kiran-sagar.kollepara@ec-nantes.fr}}
    \presentaddress{Institut de Recherche en Génie Civil et Mécanique (GeM), Ecole Centrale de Nantes, Nantes, France.} 
\fi

\def\AbsText{Low-rank model order reduction strategies for contact mechanics show limited dimensionality reduction due to linear inseparability of contact pressure field. Therefore, a dictionary based strategy is explored for creating efficient models for frictionless non-adhesive contact. A large dictionary of contact pressure trajectories is generated using a high-fidelity finite element model, while approximating the online query with a small number of dictionary entries. This is achieved by inducing sparsity in the approximation. Accuracy, computational effort and limitations of such methods are demonstrated on few numerical examples.}

\ifx\fmt\undefined
    % \newpage
    \date{}
    \maketitle
    \begin{abstract}
        \AbsText
    \end{abstract}

\else
    \abstract[Abstract]{\AbsText}


    \renewcommand\thefootnote{}
    \footnotetext{\textbf{Abbreviations:} FEM, Finite Element Method; LMM, Lagrange Multiplier Method; KKT, Karush-Kuhn-Tucker; ROM, Reduced Order Modelling, OMP, Orthogonal Matching Pursuit; nnFOCUSS, non-negative FOCUSS; SVD, Singular Value Decomposition; POD, Proper Orthogonal Decomposition}

    \renewcommand\thefootnote{\fnsymbol{footnote}}
    \setcounter{footnote}{1}
    
    \maketitle
\fi

\section{Introduction}
\label{sec:introduction}
Contact mechanics deals with multi-body mechanical problems with the possibility of the outer surfaces coming to contact. The occurrence and position of contact between bodies is generally unknown a priori, which adds non-linearity to the problem. Such problems have a wide range of applications; including tribology~\cite{Li2022}, crash simulations~\cite{Brown2000}, manufacturing~\cite{Button2013} and biomechanics~\cite{Ateshian2015} to name a few. Important contributions to numerical treatment of contact problems have been made in the past few decades with advances in algorithms as well as the interface physics~\cite{Wriggers1985,Hallquist1985,Fischer2006,Wriggers2006,Chouly2017}. The complexity of the contact model can have varying degrees, with as simple as hard contact with only compressive surface traction, to as complex as frictional and adhesive contact. Also, the search for contact pairs can be a computationally intensive task.

In the context of finite element methods (FEM), major bottlenecks include solving contact detection and associated non-linearities. Contact mechanics problems are usually variational inequalities, in other words, they are inequality constrained minimization problems. In the simplest case of non-adhesive frictionless contact, the inequality constraint is essentially a non-penetration condition. Such problems involve the computation of mainly two physical unknowns: the displacement field (primal) and a non-negative contact pressure (dual). Various methods exist to resolve the mixed problem, such as the Lagrange Multiplier method (LMM), Augmented LMM, and the mortar method, among others. In the LMM, the so-called Karush-Kuhn-Tucker (KKT) conditions impose the non-negativity of contact pressure, non-penetration of contact surfaces and finally a complimentary slackness condition which forces either the contact pressure or the gap between surfaces to be null at every position on the potential contact surface~\cite{Wriggers2006,Yastrebov2013}.

% The non-negativity of contact pressure is enforced by the so-called Karush-Kuhn-Tucker (KKT) conditions through the Lagrange Multiplier~\cite{Wriggers2006,Yastrebov2013}.

As the importance of real time mechanical simulations across industries grows rapidly~\cite{Meier2005,Tan2023}, it is important to address reduced order modelling (ROM) methods applied to contact mechanics. Usually, ROM methods typically consist of an \emph{offline} and \emph{online} stage. In the offline stage, a high-fidelity model is evaluated for different parametric instances to obtain training data, referred to as \emph{snapshots}. The reduced model is then built using the offline data and the governing equations. In the online phase, the user can evaluate the reduced model very efficiently, permitting exploration of the parametric space.

In an LMM-based approach to reduced modelling of contacts, two reduced basis are sought. This can be seen in one of the earliest contributions from~\citet{Haasdonk2012}, where the reduced contact pressure space is defined as the span of contact pressure snapshots. Compression of contact pressure data was explored by~\citet{Balajewicz2016} and \cite{Niakh_mCPG_2023,Benaceur2020} using Non-Negative Matrix Factorization and greedy algorithms with a non-negativity preserving projection, respectively. A reduced integration domain (RID) based approach is presented in~\citet{Fauque2018,LeBerre2022}, where the surface integrals and contact pressure degrees of freedom (dofs) are evaluated only on the surface of the RID in a reduced setup, followed by recovery of contact pressure through a non-negative least-squares interpolation from a training set.

Some works have bypassed the creation of dual reduced bases. For example,~\citet{Bader2016,Niakh_Nitsche_2023} use the penalty approach instead of LMM which eliminates the dual variable at the cost of allowing small penetration values. On the other hand,~\cite{Manvelyan2021} puts forward the argument that reducing the dual dofs is unnecessary as the number of surface dofs is much less than the number of internal dofs. However, such assumptions may not always hold true in examples such as membrane-like or shell-like geometries that experience contact~\cite{Haasdonk2012,Mulye2020}.

Most works mentioned above assume the presence of a low-dimensional contact pressure space for building the reduced model. However, as shown recently in our recent work~\citet{Kollepara2022}, contact problems have a peculiar feature that renders this assumption untrue in many cases. This feature is the contact zone, that exhibits a local nature. The contact zone is a function of the state of the system, and hence, may strongly depend on the loading, geometry, and other physical parameters. 

In general, for variational inequalities, only some of the constraints are usually active. An inequality constraint is referred to as active only if the current state of the problem satisfies the equality. The so-called complementarity slack KKT condition implies that the Lagrange multipliers assume a non-zero value only where constraints are active. In other words, contact pressure is non-zero only at parts of the surface where contact is established. The strong dependence of the active zone on the state of the system leads to the so-called \textit{linear inseparability} of contact pressure. In such a situation, contact pressure field cannot be approximated well using linear combinations of other snapshots. Thus, the snapshots can neither be compressed nor do they \textit{generalise}\footnote{Generalisation ability is a measure of model error outside the training set~\cite{davies2002}} well~\cite{Kollepara2022}. This has two implications from an ROM point of view:
\begin{itemize}
    \item It becomes necessary to generate a large training set of contact pressure snapshots.
    \item The low-rank contact pressure basis will have approximately the same dimension as the high-fidelity solution space; thus, offering negligible reduction.
\end{itemize}

\begin{figure}[!htb]
    \centering
    \begin{subfigure}[t]{0.49\linewidth}
        \includegraphics{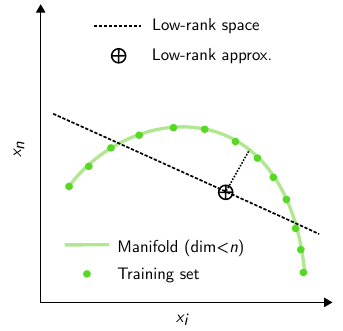}
        \caption{Low-rank approximation}
        \label{fig:manifold_approx_LR}
    \end{subfigure}
    \begin{subfigure}[t]{0.49\linewidth}
        \includegraphics{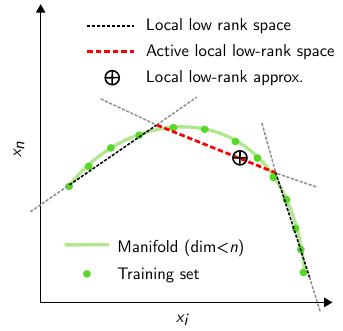}
        \caption{Localized low-rank approximation}
        \label{fig:manifold_approx_local_LR}
    \end{subfigure}\vspace{0.5em}
    \begin{subfigure}[t]{0.45\linewidth}
        \includegraphics{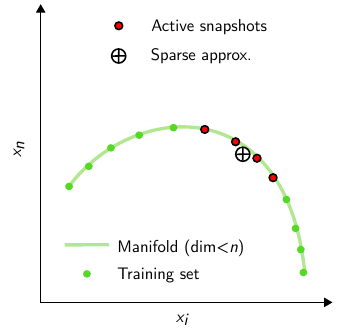}
        \caption{Sparse approximation}
        \label{fig:manifold_approx_sparse}
    \end{subfigure}
    \caption{Illustration of low-rank, localized low-rank and dictionary-based approximation. Snapshots (green circles) in an $n$-dimensional system lie on a manifold (green-curve) of dimension lower than $n$. Low-rank subspaces have a dimension $\ll n$. The approximation to the queried solution ($\oplus$) and its closeness to the manifold can be visualized for the three methods.}
    \label{fig:manifold_approx}
\end{figure}

The first problem cannot be avoided, as a large database of snapshots is required for a satisfactory interpolation of the non-linear manifold. However, it is possible to mitigate the second one, the ideas for which are discussed using the visual aid of \cref{fig:manifold_approx}. The figure illustrates the high-dimensional solution space, inside which snapshots lie in a low-dimensional solution manifold.  If the solution manifold is significantly non-linear, a low-rank approximation will generally not be satisfactory (\cref{fig:manifold_approx_LR}). To capture the manifold, a better strategy would be to generate an approximation based on snapshots nearby the queried point. This can be done in two ways:
\begin{enumerate}
    \item Multiple localized solution subspaces will better approximate the curvature of the manifold~\cite{Amsallem2012}, as depicted in \cref{fig:manifold_approx_local_LR}. However, to construct localized subspaces, a clustering technique is needed to group snapshots into sets based on distance.
    \item Sparse regression methods can be used to select the best snapshots from the training set to render a close approximation (\cref{fig:manifold_approx_sparse}). In other words, the sparse algorithm is allowed to select the locally significant snapshots, without the need for clustering. Indeed, this approach is equivalent to finding the best ``cluster'' for the queried point on-the-go. Thus, there is no need of ``discretizing'' the manifold into local subspaces. Also, the reconstruction quality in regions near intersection of two local subspaces need not be compromised.
    \item Another possibility is to use non-linear approximation techniques~\cite{Amsallem2009}, which is outside the scope of this article (and will be explored in a future article).
\end{enumerate}

In this article, we propose to use sparse approximations from the domain of dictionary learning~\cite{Aharon2006,Rubinstein2014}. Dictionary-based sparse approximations have already been applied to various numerical models including transport, acoustics~\cite{Balabanov2021} and biomedical~\cite{Lauzeral2019} applications. In contrast to a reduced basis, an over-complete dictionary contains a large number of snapshots, some of which would essentially contain redundant information. For an approximation using linear combinations of the dictionary elements, multiple solutions would likely exist. The idea is to induce sparsity in the coefficients of the online problem, encouraging selection of a few snapshots that are the best candidates for the problem. The reader should beware that the underlying linear inseparability is not resolved by this approach. While the use of large dictionaries provide a better resolution of the non-linear manifold, the primary focus of this work is overcoming the added complexity resulting from use of large dictionaries.

Our work has some similarities as well as differences with earlier approaches. While~\citet{Haasdonk2012} uses snapshots to define solution subspaces, sparse regularization is not imposed for the dual reduced variable. Moreover, our approach does not attempt to compute a dual reduced basis like the works of~\citet{Balajewicz2016,Benaceur2020,Niakh_mCPG_2023}. Both~\citet{Haasdonk2012,Fauque2018} tackle problems with node-to-node interaction, whereas our work includes problems with node-to-segment interactions and large ``slip'' between contact surfaces.

The article is structured as follows: first a very brief introduction to the discretized contact problem is given, followed by a brief introduction to sparse regression methods. Then, two proposed formulations of dictionary-based sparse methods to contact mechanics are introduced with numerical examples. Finally, some observations on the two methods are discussed and analyzed, prior to concluding the article.

\section{Model Problem}
\begin{figure}
    \centering
    \includegraphics[scale=0.7]{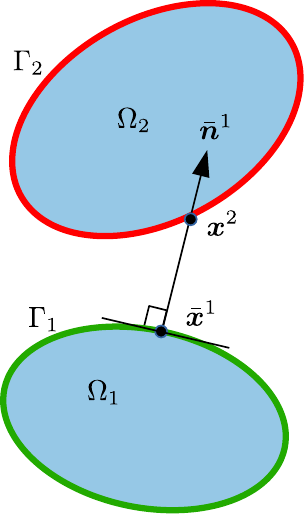}
    \caption{Kinematic description of two body contact problems with possibility of contact}
	\label{fig:contact-schematic}
\end{figure}

Typical contact problems can be framed as inequality constrained minimization problems. Consider two elastic bodies with frictionless and non-adhesive surfaces, that might come into contact in a deformed state (\cref{fig:contact-schematic}). The contact problem can be framed as a minimization of the internal elastic energy while respecting the non-penetration constraints. To remain concise, the discretized version of this problem is given here. For a more detailed formulation, the reader may refer to the preceding publication~\cite{Kollepara2022} as well as the literature~\cite{Wriggers2006,Wriggers1985,Hallquist1985}.

\begin{gather}
    \begin{aligned}
        \argmin_{\vh} \frac{1}{2} {\vh}^T \mathbf{K}(\mu) \vh - {\vh}^T \bm{f}(\mu) \\
        \mathbf{C}(\mu,\vh) \vh - \bm{g}(\mu,\vh) \leq \bm{0} 
    \end{aligned}
    \label{eq:inequality_discrete}
\end{gather}

\noindent where 
\begin{itemize}
    \item $\vh = \begin{bmatrix} \bm{v}^h_1 \\ \bm{v}^h_2 \end{bmatrix}$ is the stacked nodal displacement vector
    \item $\mathbf{K}(\mu) = \begin{bmatrix} \mathbf{K}_1(\mu) & 0 \\ 0 & \mathbf{K}_2(\mu) \\ \end{bmatrix}$ is the block elastic stiffness matrix
    \item $\bm{f}(\mu)  = \begin{bmatrix} \bm{f}_1(\mu) \\ \bm{f}_2(\mu) \end{bmatrix}$ is the stacked nodal loading vector
    \item $\bm{g}(\mu,\vh)$ is the node-to-segment distance between the surfaces of two bodies in reference configuration \footnote{This quantity depends on $\vh$ because contact pairs are established in deformed configuration. See~\cite{Kollepara2022} for details.}
    \item $\mathbf{C}(\mu, \vh) \vh$ gives the change in distance between the surfaces in the deformed configuration.
\end{itemize}
with the subscripts $1$ and $2$ referring to the two bodies.

\noindent Solution to such an inequality constrained problem is given by the following Karush-Kuhn-Tucker (KKT) conditions, where an additional variable i.e.\ the Lagrange Multiplier $\lambdah$ is introduced:

\begin{subequations}
	\begin{align}
        \mathbf{K}(\mu) \uh - \bm{f}(\mu) + \mathbf{C}^T(\mu,\uh) \lambdah &= \bm{0}  \label{eq:kkt_force_balance} \\
        \lambdah &\geq \bm{0} \label{eq:kkt_nonnegativity}  \\
        \mathbf{C}(\mu,\uh) \uh - \bm{g}(\mu,\uh) &\leq \bm{0} \label{eq:kkt_nonpenetration} \\
        {\lambdah}^T [\mathbf{C}(\mu,\uh) \uh - \bm{g}(\mu,\uh)] &= 0 \label{eq:kkt_comp_slack}
	\end{align}  \label{eq:kkt_discrete}
\end{subequations}

\noindent In this case, the Lagrange Multiplier $\lambdah$ corresponds to the contact pressure between the two surfaces.

\section{Dictionary learning with sparse methods: a brief overview}\label{sec:sparse_methods}
Before getting started with the contact problem, sparse regression techniques and corresponding literature are briefly introduced. This is followed by a brief discussion on their application to partial differential equations (PDEs).

\subsection{Sparse regression methods}
Sparse regression techniques have been under development  in recent decades across various fields, especially in signal processing and statistics\cite{Simoncelli1992,Candes2002}. The main idea behind these techniques is to approximate a signal with the few most suitable elements of a precomputed dictionary of signals~\cite{Elad2010}, or to develop the simplest predictive model~\cite{Brunton2016}. A typical sparse regression can be written as the following optimization statement:
\begin{align}
    \centering
    \begin{aligned}
        \min & \norm{\vphantom{\widehat{}} \bm{\alpha}}{p} \\
        \text{s.t.} & \norm{\bm{x} - \mathbf{D} \bm{\alpha}}{}  < \varepsilon
    \end{aligned}
    \label{eq:parsimony_idea}
\end{align}
where $\bm{x}$ is the signal being approximated, $\mathbf{D}$ could either be a dictionary of signals or a basis that possibly admits sparse representation, $\bm{\alpha}$ is a vector of coefficients, $\norm{\cdot}{p}$ is a sparsity inducing norm, and $\varepsilon$ is the tolerance on the approximation error. It is also possible to write other forms of sparse regression problems where the approximation error is minimized and $p$-norm term is constrained to stay under a specific tolerance of sparsity or the minimization of a penalty form with a weighted sum of the two terms. A primer on sparse regression methods can be found in~\citet{Mairal2014}.

It is common knowledge in the communities using sparse regression methods that the $\ell_0$-norm is a measure of cardinality, hence~\eqref{eq:parsimony_idea} with $p=0$  is the natural choice in a theoretical sense. However, this problem is NP-hard in general, and therefore, it is more common to use $\ell_p$-norm with $p\in ]0,1]$. Examples of such methods include the LASSO~\cite{Tibshirani1996} and Dantzig Selector~\cite{Candes2007}, which use $p=1$. Other possibilities include the use of greedy algorithms, such as Orthogonal Matching Pursuit (OMP)~\cite{Pati1993} and weight $\ell_2$-norm minimization such as FOCUSS~\cite{Gorodnitsky1997}. OMP and FOCUSS methods, which will be referred to in upcoming sections, are briefly detailed in \cref{sec:sparse_methods_list}.

\subsection{Numerical modelling using sparse methods} \label{sec:dictionary_resolution}
Sparse regression methods, though originally developed for application on pure data approximation or data compression problems, have seen recent interest in applications towards physics problems~\cite{Wang2011,Schaeffer2013,Brunton2016}. The methods that have been discussed previously are also meant for data approximation. At this point, the application of such methods to unconstrained PDE is discussed, which can be formulated as a residual minimization problem involving dictionary-based approximation and sparse regression. Consider a discretized problem with the state vector $\bm{v}$ and residual $\bm{r}(\bm{v})$. Given a dictionary $\mathbf{D}$, the basic idea is find an approximation $\bm{v} \approx \mathbf{D} \widehat{\bm{v}}$ that minimizes a sparsity inducing $\ell_p$-norm of $\widehat{\bm{v}}$ but also keeping the residual $\bm{r}(\mathbf{D} \widehat{\bm{v}})$ is below a certain tolerance, as given in~\eqref{eq:sparse_res}.

\begin{gather}
    \begin{aligned}
        \min  & \norm{\widehat{\bm{v}}}{p} \\
        \text{s.t. } &\norm{\bm{r}(\mathbf{D} \widehat{\bm{v}})}{} < \varepsilon
    \end{aligned}
    \label{eq:sparse_res}
\end{gather}

Dictionary-based sparse approximation for PDEs has been discussed in~\cite{Balabanov2021}, along with the application of random sketching to residual evaluation which provides an efficient way of solving the minimization problem. This is similar to the approach of randomized SVD proposed in~\cite{Halko2011}, where random combinations of rows/columns are used. The idea is to project the residual on a low-dimensional subspace in $\mathtt{colsp}(\mathbf{D})$, where $\mathtt{colsp}$ denotes the column-space. The matrix $\mathbf{B}$, containing a basis corresponding to a random linear combination of dictionary columns, is used for this projection in~\eqref{eq:sparse_res_randomized}. This projection essentially reduces the complexity of enforcing the residual constraint. Given a dictionary $\mathbf{D}$ of size $m\times n$, a smaller matrix that has information from all columns of $\mathbf{D}$ can be generated by multiplying a matrix $\mathbf{R}$ of size $n\times L$ with $L \ll \min{(m, n)}$. After orthogonalization, the output matrix $\mathbf{B}$ can be used for projecting the residual $\bm{r}$.
\begin{subequations}
    \begin{align}
        \min & \norm{\widehat{\bm{v}}}{p} \\
        \text{s.t. }&\norm{\mathbf{B}^T \bm{r}(\mathbf{D} \widehat{\bm{v}})}{} < \varepsilon \\
        \text{where } &\mathbf{B} = \mathtt{orth}(\mathbf{D} \mathbf{R}) \ , \ \ \mathbf{R} \in \mathcal{U}_{[0,1]}^{n \times L} \label{eq:randomized_orth_matrix}
    \end{align}
    \label{eq:sparse_res_randomized}
\end{subequations}
where $\mathtt{orth}(\cdot)$ denotes an orthogonalizing operation and $\mathcal{U}_{[0,1]}^{n \times L}$ denotes a random matrix of size $n \times L$ whose entries are drawn on the interval $[0,1]$ using a uniform distribution. The number $L$ can be tuned to not only control the complexity of the minimization problem but also to crudely constrain the sparsity of the solution. For $L \ll n$, $\mathtt{colsp}(\mathbf{B}) \subset \mathtt{colsp}(\mathbf{D}$) is low-dimensional and hence, the projection $\mathbf{B}^T \bm{r}$ is a ``weaker'' evaluation of $\bm{r}$, which allows for a sparser $\widehat{\bm{v}}$ at the price of permitting residuals orthogonal to $\mathbf{B}$.

Another related approach in which snapshots are directly used for reducing the model without any explicit compression is the so-called CUR decomposition~\cite{LeGuennec2018}. In this method, not only columns but also rows from the snapshot matrix are used for developing a parametric regression. It works by algorithmically selecting a column and row subset from the snapshot matrix, and using them for developing a regression model. A contrasting feature of the CUR approach is the use of regression without invoking the partial differential equation or its weak form, resulting in a non-intrusive method.

In the proposed dictionary-based approach for contact mechanics, unlike the low-rank approach~\cite{Balajewicz2016,Benaceur2020}, no compression of the contact pressure snapshots is applied. The dual solution can simply be restricted to the column-cone\footnote{column-cone can be thought of as the subset of the column-space restricted to non-negative coefficients.} of the dictionary $\mathbf{D}$, which might be reasonably rich to approximate the true contact pressure subcone. In the online phase, only a few snapshots will be selected to estimate the contact pressure field. This is justified because not all information in the dictionary is necessary for the reconstruction of a particular instance and hence, the dual reduced dofs will admit a sparse solution. Therefore, enough motivation exists to use sparsity-enforcing methods to choose a few columns of a large dictionary to approximate the dual solution.

\section{Dictionary-based sparse methods for contact problems}\label{sec:sparse_dict}
This section introduces a dictionary-based sparse approximation algorithm for inequality-constrained problems. For the contact pressure dictionaries, the term ``over-complete'' is used loosely. It is simply used to imply that the dictionary has numerous snapshots, possibly larger than the number of contact dofs in the high-fidelity problem.

\begin{remark}
    The authors investigated another idea to solve the dictionary-based sparse problem. It involved strong assumptions about the convexity of the feasible region defined by the non-penetration constraint. Therefore, the method is not included in the main article because of its limited scope. However, it is briefly explained in the \cref{sec:convex_hull} with an illustration.
\end{remark}

% \subsection{A greedy active-set method for dual dictionary element selection}\label{sec:greedy}

A greedy method inspired by the Orthogonal Matching Pursuit (OMP) (see \cref{sec:sparse_methods_list}) and the active set method is devised, referred to as \emph{greedy active-set} method. The primal and dual reduced spaces are defined using a reduced basis and a dictionary, as detailed below : 
\begin{itemize}
    \item Primal reduced basis $\PrimalRB$ is defined using the left singular vectors of primal snapshots, i.e.\ $\PrimalRB \gets \text{svd}(\PrimalDict, \delta)$ where $\delta$ is the truncation tolerance. The first singular vectors that contribute to $1-\delta$ fraction of total energy (computed using squared singular values) are kept.
    \item The dual dictionary $\DualDict$ contains all the snapshots from the training phase. 
\end{itemize}
    The above choice is justified since the displacement field is known to be relatively separable compared to the contact pressure\cite{Kollepara2022}. The unknowns of the online problem are the reduced displacement dofs $\widehat{\bm{u}} $ and the sparse contact pressure dofs $\widehat{\bm{\lambda}}_{\mathcal{I}}$. Using the defined notations, a projection-based, sparse problem can be stated as:

\begin{subequations}
	\begin{align}
		&\min \norm{\widehat{\bm{\lambda}} }{0} \\ %  \label{eq:sparse_greedy_minimize}\\
		& \text{s.t. }  \begin{bmatrix}
			\mathbf{\PrimalRB}^T \mathbf{K} \mathbf{\PrimalRB}  & \mathbf{\PrimalRB}^T \mathbf{C}^T \DualDictActive \\
			\DualDictActive^T \mathbf{C} \mathbf{\PrimalRB} & \mathbf{0}
		\end{bmatrix} 
		\begin{bmatrix}
			\widehat{\bm{u}} \\
			\widehat{\bm{\lambda}}_{\mathcal{I}}
		\end{bmatrix} =  
        \begin{bmatrix}
			\mathbf{\PrimalRB}^T \bm{f} \\
			\DualDictActive^T \bm{g}
		\end{bmatrix} \label{eq:sparse_greedy_system} \\
		&  \text{and } \DualDict^T  \left ( \mathbf{C} \mathbf{\PrimalRB} \widehat{\bm{u}} - \bm{g} \right ) \leq \tau  \label{eq:sparse_greedy_kkt_nonpen}\\
        &   \phantom{\text{and }} \widehat{\bm{\lambda}} \geq 0 \label{eq:sparse_greedy_kkt_nnpr}\\
		& \text{where $\mathcal{I}$ indicates the active set and }  \DualDictActive = \DualDict[:,\mathcal{I}] \nonumber
	\end{align}
\label{eq:sparse_greedy}
\end{subequations}

\noindent\eqref{eq:sparse_greedy_system},~\eqref{eq:sparse_greedy_kkt_nonpen}~and~\eqref{eq:sparse_greedy_kkt_nnpr} are based on the same principle as the KKT conditions~\eqref{eq:kkt_discrete}, but applied over solution spaces defined by the reduced basis $\PrimalRB$ and the dictionary $\DualDict$ (rather than a low-rank dual basis). Comparing~\eqref{eq:sparse_greedy} with the generic form of sparse regression statement~\eqref{eq:sparse_res}, the balance equation is imposed as an equality constraint, instead of an inequality constraint where the residual norm is limited to a certain tolerance. Imposing limits on the residual norm is usually more relevant when using random sketching, like in~\eqref{eq:sparse_res_randomized} when the system is under-determined. Also, it is not straightforward to use standard sparse regression algorithms to solve inequality constraints~\eqref{eq:sparse_greedy_kkt_nonpen}~and~\eqref{eq:sparse_greedy_kkt_nnpr}.

%Equality constraint can be used as long as the system is determined, which is currently the case in~\eqref{eq:sparse_greedy_system}.

The truncated primal basis $\PrimalRB$ cannot contain the entire information learnt from the training set. The solutions generated using the resulting subspace tend to have some \textit{false} penetrations, which tend to drive the selection of some additional but irrelevant contact pressure snapshots from the dictionary $\DualDict$ (see \cref{sec:tau_delta}). To avoid this issue, a small penetration $\tau$ is permitted in \eqref{eq:sparse_greedy_kkt_nonpen} to avoid this spurious selection of snapshots. Since the unwanted spurious effects were observed to increase with truncation, the permitted penetration $\tau$ must be linked to the tolerance $\delta$. Unless specified, the $\tau := \delta$ is used in the rest of the article.

\def\remarktext{
    As a consequence of~\eqref{eq:sparse_greedy_kkt_nonpen}, the non-penetration condition is not satisfied exactly by the sparse model. Thus, one might question the interest of using of an LMM-based model over a penalty-based unconstrained model that is simpler to reduce. While this is a valid argument, the community remains interested in reduction of LMM-based models~\cite{Haasdonk2012,Balajewicz2016,Fauque2018,Benaceur2020,LeBerre2022}, possibly for reasons such as avoiding ill-conditioning problems of penalty methods.
}

\ifx\fmt\undefined
    \begin{remark}
        \remarktext
    \end{remark}
\else
    \begin{boxwithhead}
    {Remark} 
    {\noindent \remarktext}
    \end{boxwithhead}
\fi

The primary notion in greedy algorithms is to enrich the approximation by a rank-1 term in each step. For standard sparse regression problems like~\eqref{eq:sparse_res}, the OMP algorithm starts with an empty coefficient vector and adds a single non-zero entry in one step. The position of the new non-zero entry corresponds to the vector in dictionary with the highest correlation to the residual. Using the same idea, the greedy active-set algorithm devised here performs enrichment by correlating the contact pressure snapshots with the violation of inequality constraint~\eqref{eq:sparse_greedy_kkt_nonpen}. More precisely, the algorithm starts with a zero $\widehat{\bm{\lambda}}$ and in each enrichment step, the dictionary column with the highest correlation to the violations of the non-penetration condition is added in the greedy enrichment. Thus, enrichment index of the dictionary $p_{\texttt{add}}$ can be computed as:

\begin{align}
    p_{\texttt{add}} &:= \argmax_p \DualDict[:,p]^T  \left ( \mathbf{C} \PrimalRB \widehat{\bm{u}}^{k} - \bm{g} \right )^{+, \tau}
    \label{eq:sparse_enrich_concept}
\end{align}
 where  $
	(\bm{z})^{+,\tau}_i :=\left \lbrace\begin{matrix}
		z_i & \text{if } z_i>\tau \\
		0 & \text{otherwise}
    \end{matrix}  \right . $, is a hard-thresholding operator giving the part greater than $\tau$. \\

\noindent\eqref{eq:sparse_enrich_concept} can be computed more efficiently using the reduced operators $\widehat{\mathbf{C}}(\mu,  \widehat{\bm{u}}^{k-1})$ and $\widehat{\bm{g}}(\mu,  \widehat{\bm{u}}^{k-1})$:
\begin{align}
    p_{\texttt{add}} &:= \argmax_p  \left ( \widehat{\mathbf{C}} \widehat{\bm{u}}^{k} - \widehat{\bm{g}} \right )^{+,\tau}[p] % \label{eq:sparse_enrich_rb}
\end{align}
where 
\begin{align*}
    \widehat{\mathbf{C}}(\mu,  \widehat{\bm{u}}^{k-1}) &= \DualDict^T \mathbf{C}(\mu,  \widehat{\bm{u}}^{k-1}) \PrimalRB  \\
    \widehat{\bm{g}}(\mu,  \widehat{\bm{u}}^{k-1}) &= \DualDict^T\bm{g}(\mu,  \widehat{\bm{u}}^{k-1})
\end{align*}

However, to satisfy the inequality constraint~\eqref{eq:sparse_greedy_kkt_nnpr} that prohibits non-negative pressures, it is necessary to not only enrich but also to eliminate terms that violate the constraints in each greedy step. This is done by eliminating the largest negative $\widehat{\bm{\lambda}}$ i.e.\ setting the largest negative entry to zero. In the absence of a negative entry, elimination is not performed. The elimination index $p_{\texttt{rem}}$ is computed as:
\begin{align}
    p_{\texttt{rem}} &= \argmin_p \left ( \widehat{\bm{\lambda}}^{k} \right )^{-} [p]
\end{align}

\noindent where $
	(\bm{z})^-_i :=\left \lbrace\begin{matrix}
		z_i & \text{if } z_i<0 \\
		0 & \text{otherwise}
	\end{matrix}  \right . $

    Although the motivation for this algorithm was based on the OMP, the greedy active-set algorithm ultimately resembles a fixed-point active-set method. However, a maximum of one dual dictionary column is activated in a given iteration, thereby maintaining a small active set in intermediate iterations. The algorithm is given in \cref{alg:sparse_greedy}, where Steps~\ref{step:enrich_start} to~\ref{step:enrich_end} contain the enrichment and the elimination process. In each iteration, only one of the two operations is performed. Due to the architecture of the algorithm, the complexity of operations does not scale with the dictionary size, except enrichment/elimination Steps~\ref{step:enrich_start} to~\ref{step:enrich_end}. The small sizes of active set $\mathcal{I}$, which is directly related to the sparsity of $\widehat{\bm{\lambda}}$, prevent the size of mixed system in Step~\ref{step:invert} from growing significantly during the iterative process.

    Note that at any given greedy step, $p_{\texttt{rem}}$ will never be same as $p_{\texttt{add}}$. For any active dual dof indexed $\widetilde{p}\in \mathcal{I}$, \eqref{eq:sparse_greedy_system} implies $\widehat{\mathbf{C}}[\hspace{0.15em}\widetilde{p}\hspace{0.15em}] \widehat{\bm{u}}^k - \widehat{\bm{g}}[\hspace{0.15em}\widetilde{p}\hspace{0.15em}] = 0 \implies p_{\texttt{add}} \neq \widetilde{p} \implies p_{\texttt{add}} \neq p_{\texttt{rem}}$, since $p_{\texttt{rem}} \in \mathcal{I}$. Moreover, to avoid any possibility of an endless loop between same $p_{\texttt{add}}$ and $p_{\texttt{rem}}$, an upper limit to number of iterations $k_{max}$ is added.

\begin{algorithm}[!htb]
    \caption{Greedy active-set algorithm} \label{alg:sparse_greedy}
    \begin{algorithmic}[1]
        \State{Input: Queried value of parameter $\mu$}
        \State{Given: Primal basis $\PrimalRB$ and dual Dictionary $\DualDict$}
        \Statex{\hspace{3em} Reduced operators $\mathbf{\PrimalRB}^T \mathbf{K} \mathbf{\PrimalRB} \text{ and } \mathbf{\PrimalRB}^T \bm{f}$} \Comment{can be built offline}
        \State{Initialize: $k = 0, \ \ \mathcal{I} = \emptyset$}
        \While{$\widehat{\bm{u}}$ and $\mathcal{I}$ not converge AND $k < k_{\text{max}}$} \label{step:loop_start}
            \State{Build reduced constraint operators $\widehat{\mathbf{C}}(\mu,  \widehat{\bm{u}}^{k-1})$ and $\widehat{\bm{g}}(\mu,  \widehat{\bm{u}}^{k-1})$} \label{step:build_operators}
            \State{Solve the linear system~\eqref{eq:sparse_greedy_system}} \label{step:invert}
            \State{Set $\widehat{\bm{\lambda}}^{k}[\mathcal{I}] \gets \widehat{\bm{\lambda}}_{\mathcal{I}}$ and $\widehat{\bm{\lambda}}^{k}[\mathcal{I}^{\mathsf{c}}] \gets \bm{0}$} \Comment{${\mathcal{I}}^{\mathsf{c}}$ is the complementary set of $\mathcal{I}$} 
            \Statex{}
            \Statex{\hspace{1em} Enrich/Eliminate decision :}
            \If{$ \widehat{\bm{\lambda}}^{k}  \geq \bm{0}$ :} \label{step:enrich_start} \Comment{All multipliers are non-negative, so enrich} \vspace{0.8em}
            \Statex{\centerline{$ \begin{array}{rcl}
                    p_{\texttt{add}} &=& \argmax_p  \left ( \widehat{\mathbf{C}} \PrimalRB \widehat{\bm{u}}^{k} - \widehat{\bm{g}} \right )^{+,\tau}  \\
                    \mathcal{I} &\gets&  \mathcal{I} \cup \{ p_{\texttt{add}}  \}   
            \end{array} $ }} 
            \Else: \Comment{Eliminate negative multiplier}
            \Statex{\centerline{$ \begin{array}{rcl}
                    p_{\texttt{rem}} &=& \argmin_p \left ( \widehat{\bm{\lambda}}^{k} \right )^{-} [p] \\
                    \mathcal{I} &\gets&  \mathcal{I} \setminus  \{ p_{\texttt{rem}} \}
            \end{array} $ }}
            \EndIf{} \label{step:enrich_end}
            \State{$k \gets k+1$}
        \EndWhile{} \label{step:loop_end}
        \State{Reconstruct $\bm{u} = \PrimalRB \widehat{\bm{u}}$ and $\bm{\lambda} = \DualDict \widehat{\bm{\lambda}}$}
        \State{Output: $\bm{u}, \bm{\lambda}$}
    \end{algorithmic}
\end{algorithm}

\subsection{Application to the Hertz problem}\label{sec:greedy_hertz} 
The Hertz contact problem between two half cylinders (\cref{fig:hertz}) is considered, with $R_1 = R_2 = 1$. A finite element model with 466 quad elements in each half-cylinder is used as a high-fidelity model. Bilinear shape functions describe the displacement field, while node-centered piece-wise constant shape functions at the boundary describe the contact pressure field. Surface nodes also serve as Gauss quadrature points for boundary integration. The imposed displacement on the top cylinder $d \in (0,0.3)$ is considered as the parameter. {Each half-cylinder has $513$ nodes and $466$ elements of which $78$ elements lie of the potential contact surface, that is, the semi-circular edges. Computational time for each snapshot is around $0.4$-$1.0$s.\\

\begin{figure}[!htb]
	\centering
    \begin{subfigure}[t]{0.44\linewidth}
        \includegraphics[width=1\linewidth]{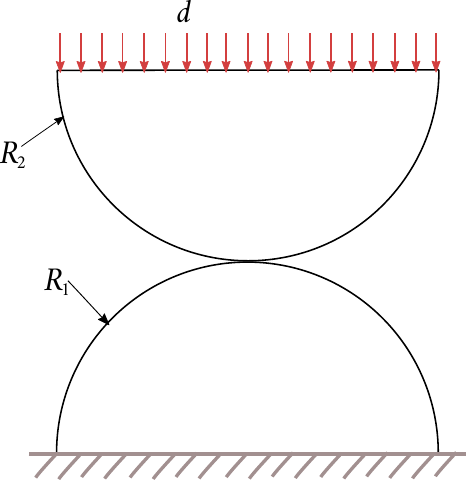}
        \caption{Hertz problem}\label{fig:hertz_half_cyl}
    \end{subfigure}
    \begin{subfigure}[t]{0.55\linewidth}
        \includegraphics[width=1\linewidth]{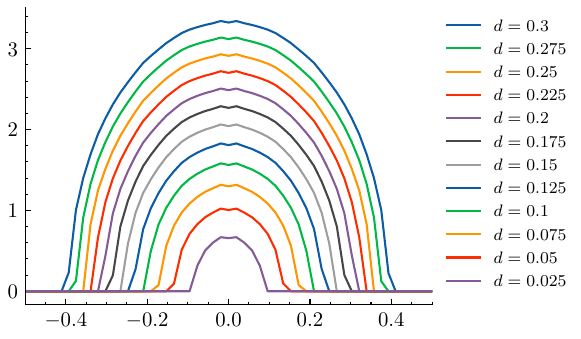}
        \caption{Contact pressure snapshots}\label{fig:hertz_snaps}
    \end{subfigure}
    \caption{Hertz problem: two half cylinders loaded against each other and the resulting contact pressure snapshots. Displacement $d$, imposed on the top cylinder, is treated as a parameter in the reduced model.}
    \label{fig:hertz}
\end{figure}

\noindent Details of the offline phase : 
\begin{itemize}
    \item Training set: Primal and dual snapshots were computed for four training sets of sizes 12, 30, 60, and 120 in the offline stage. The training points $\mathcal{P}_\texttt{tr}$ are distributed uniformly in parametric space $\mathcal{P} = (0,0.3]$.
    \item Validation Set: $\mathcal{P}_\texttt{val} \subset \mathcal{P}$ consisting 119 points, which are also the mid-points of the fourth (and largest) training set points, is used to study reconstruction errors.
\end{itemize}

The primal basis is computed using SVD with truncation parameter $\delta$, while the dual dictionary is composed of contact pressure snapshots. A spy plot of the dual dictionary is shown in \cref{fig:hertz_dictionary}.

\begin{figure}[!htb]
	\centering
    \includegraphics[width=0.5\linewidth]{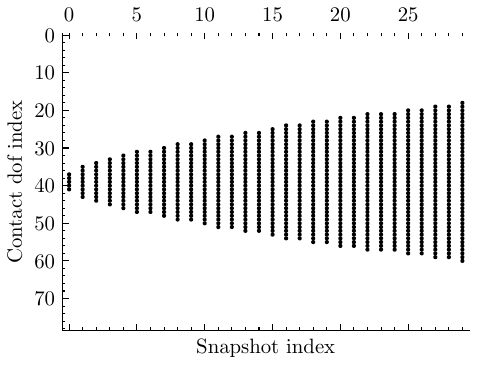}
    \caption{A spy pattern of dual dictionary with 30 elements for the Hertz problem. The $79$ points on the $y$-axis are the dofs on the potential contact surface. The dictionary columns are arranged in the increasing order of loading parameter $d$.}
	\label{fig:hertz_dictionary}
\end{figure}

Using the greedy active-set method, two instances of contact pressure reconstructions are shown in \cref{fig:hertz_reconstruction_examples}. The dictionary snapshots selected by the algorithm for reconstruction are also displayed. As expected, only a few dictionary elements, two in this case, are chosen. These two selected snapshots match the contact position required for the specific instance, indicating that the enrichment and elimination process is effectively selecting the correct snapshots for reconstruction. 

\begin{figure}[!htb]
	\centering
	\begin{subfigure}[t]{0.48\linewidth}
        \centering
        \includegraphics[width=1\linewidth]{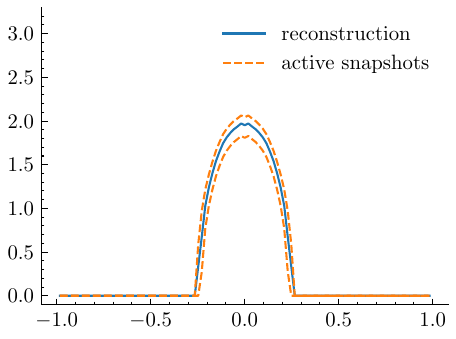}
        % \caption{$d = 0.14, \varepsilon = 1.6\times 10^{-2}$}
        \caption{{\tabular[t]{@{}l@{}}$d = 0.14$\\ Recon.\ error $1.6\times 10^{-2}$\endtabular}}
	\end{subfigure}
	\begin{subfigure}[t]{0.48\linewidth}
        \centering
        \includegraphics[width=1\linewidth]{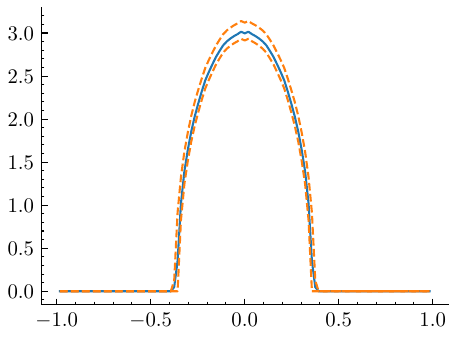}
        % \caption{$d = 0.26, \varepsilon = 2.5\times 10^{-3}$}
        \caption{{\tabular[t]{@{}l@{}}$d = 0.26$\\ Recon.\ error $2.5\times 10^{-3}$\endtabular}}
	\end{subfigure}
    \caption{Greedy active-set reconstruction of contact pressure for certain parametric instances in the validation set of the Hertz problem. Elements of the dictionary chosen by the algorithm are shown and the reconstruction errors are given. Dictionary of size 12 is used in these examples.}\label{fig:hertz_reconstruction_examples}
\end{figure}

Mean reconstruction errors over the validation set improve by more than an order of magnitude as dictionary size is increased from $12$ to $120$ for the case $\delta=10^{-10}$ (\cref{fig:hertz_mean_error_valset}). However, such gains with over-complete dictionaries diminish for the case $\delta=10^{-6}$, probably because the primal truncation starts to become the dominant source of error, rather than the lack of contact pressure separability.

Mean computational time on the validation set (\cref{fig:hertz_computation_time}) initially increases with the dictionary size, but later stabilises. However, it is evident that the computational time is primarily influenced primarily by the number of iterations (\cref{fig:hertz_num_iters}). This observation is further supported by \cref{fig:hertz_time_per_iter}, where the computational time per iteration remains relatively constant w.r.t.\ the dictionary size. The small increases in computational time per iteration with dictionary size can be explained due to the increase in the dimension of primal basis $\PrimalRB$, especially in the case of $\delta=10^{-10}$ where the primal rank is strongly dependent on the training set size.

\begin{figure}[!htpb]
	\centering
    \begin{subfigure}[b]{1.0\linewidth}
        \centering
        \includegraphics[width=0.98\linewidth]{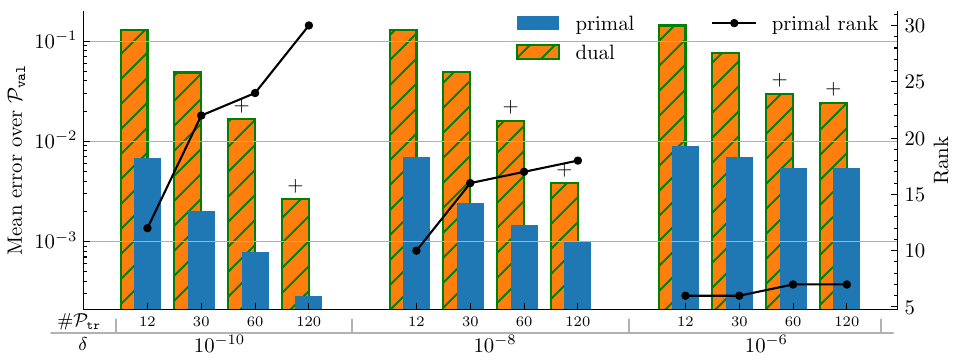}
        \caption{Mean relative errors (bar plot, left $y$-axis) and primal rank (line plot, right $y$-axis)\label{fig:hertz_mean_error_valset}. $+$ labels indicate cases where $\# \mathcal{P}_{\text{tr}} > \#$ FOM contact dofs.}
    \end{subfigure}
    \begin{subfigure}[t]{0.48\linewidth}
        \centering
        \includegraphics[width=1.0\linewidth]{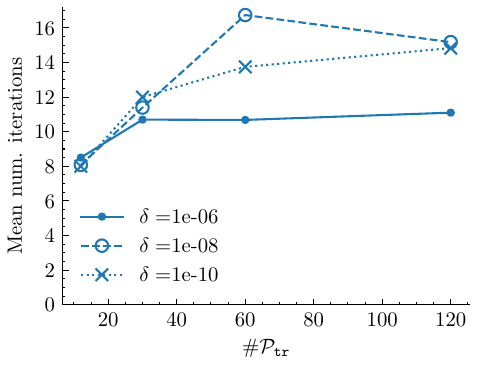}
        \caption{Mean number of iterations\label{fig:hertz_num_iters}} 
    \end{subfigure}
    \begin{subfigure}[t]{0.48\linewidth}
        \centering
        \includegraphics[width=1.0\linewidth]{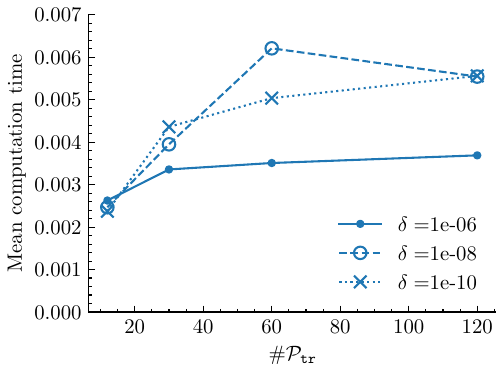}
        \caption{Mean computation time (s), (mean high-fidelity solver time $\approx 0.65$s)\label{fig:hertz_computation_time}} 
    \end{subfigure}
    \begin{subfigure}[b]{0.48\linewidth}
        \centering
        \includegraphics[width=1.0\linewidth]{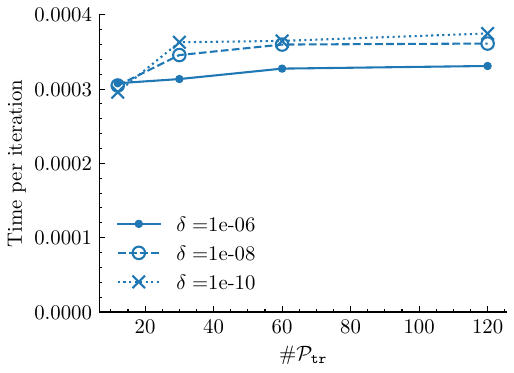}
        \caption{Mean computation time per iteration (s)\label{fig:hertz_time_per_iter}} 
    \end{subfigure}
    \caption{Evolution of \textbf{(a)} relative reconstruction errors and primal rank, \textbf{(b)} number of iterations, \textbf{(c)} total and \textbf{(d)} per iteration computation time with training set size for the Hertz problem. Mean values of these quantities over the validation set $\mathcal{P}_\texttt{val}$ are shown. Primal and dual errors are computed using $\mathcal{H}^1$ and $\mathcal{L}^2$ norm, respectively. Quantities are shown for different primal truncation tolerance  $\delta \in \{10^{-6}, \ 10^{-8}, 10^{-10}\} $. Note that computation times exclude the construction time of non-linear operators.}
	\label{fig:hertz_mean_error_and_rank}
\end{figure}

Detailed plots of reconstruction errors over the validation set $\mathcal{P}_\texttt{val}$ are shown in \cref{fig:hertz_error_valset} for primal truncation tolerance $\delta \in \{10^{-10}, 10^{-8}, 10^{-6}\}$; and each curve corresponds to one of the training sets. An unexpected observation in the case of $\delta=10^{-6}$ is that certain regions have higher dual errors for larger training sets. This happens due to the spurious penetrations discussed previously (\cref{sec:sparse_dict}) and in \cref{sec:tau_delta}. Even after relaxation of the non-penetration constraint by the value of $\tau$ (see~\eqref{eq:sparse_greedy_kkt_nonpen}), spurious selections are still possible for higher values of $\delta$. Moreover, the errors tend to be larger in general near $d=0.0$, due to higher sensitivity w.r.t.\ $d$ in this region. The cross ($\bm{\times}$) markers, indicating cases where algorithm does not converge due to oscillations between $p_{\texttt{add}}$ and $p_{\texttt{rem}}$.

\begin{figure}[!htpb]
	\centering
    \begin{subfigure}[t]{1.0\linewidth}
        \centering
        \includegraphics[width=0.48\linewidth]{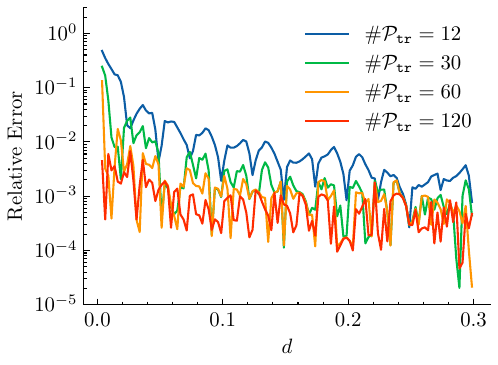}
        \quad
        \includegraphics[width=0.48\linewidth]{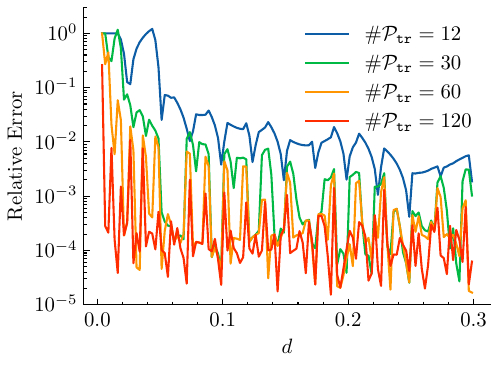}
        \caption{Primal and dual errors for $\delta = 10^{-10}$}
    \end{subfigure}
    \begin{subfigure}[t]{1.0\linewidth}
        \centering
        \includegraphics[width=0.48\linewidth]{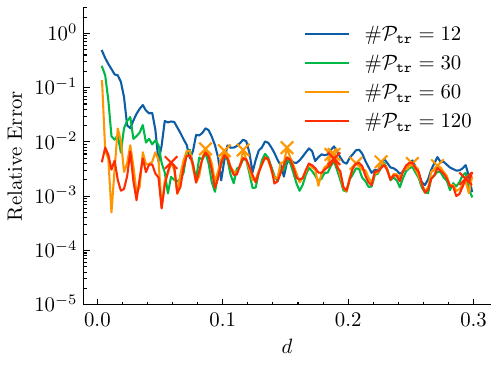}
        \quad
        \includegraphics[width=0.48\linewidth]{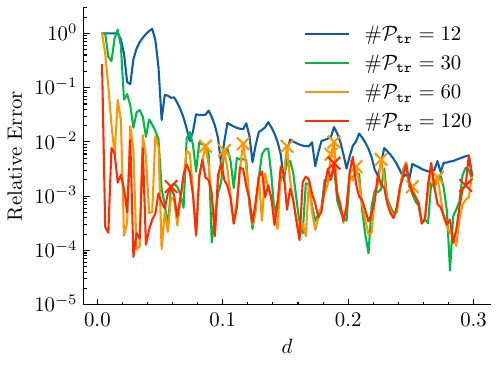}
        \caption{Primal and dual errors for $\delta = 10^{-8}$}
    \end{subfigure}
    \begin{subfigure}[t]{1.0\linewidth}
        \centering
        \includegraphics[width=0.48\linewidth]{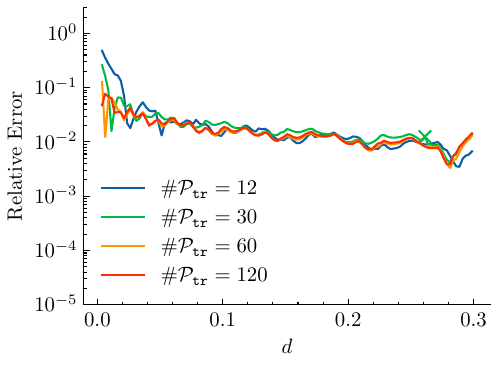}
        \quad
        \includegraphics[width=0.48\linewidth]{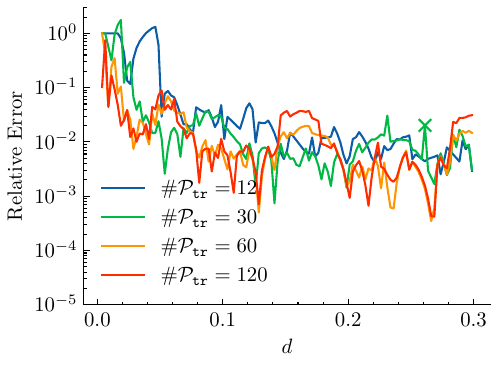}
        \caption{Primal and dual errors for $\delta = 10^{-6}$}
    \end{subfigure}
    \caption{Detailed primal and Dual reconstruction errors for the Hertz problem using the training sets 1,2 and 3. The validation set is common to all curves, consisting of 60 equispaced points in $\mathcal{P}$. Each subfigure corresponds to a different truncation tolerance $\delta$ for the primal basis. Primal and dual errors are computed using $\mathcal{H}^1$ and $\mathcal{L}^2$ norm, respectively. Cross-marks $\bm{\times}$ indicate points where the greedy active-set algorithm reached max the number of iterations ($50$) and did not converge to the defined tolerance $10^{-5}$}
	\label{fig:hertz_error_valset}
\end{figure}

The sparsity pattern of the dual reduced dofs $\widehat{\bm{\lambda}}$ for the training set of size 30 is shown in \cref{fig:hertz_dualDof_sparsity}. As the snapshots in $\DualDict$ are arranged in increasing order of $d$, most of the columns of the sparsity pattern have two close non-zeros, indicating that the nearest two snapshots to the targeted reconstruction were chosen by the algorithm. Spurious selections are seen in some rows, but the values of corresponding coefficients are quite small and do not have any significant impact on the solution.

\begin{figure}[!htb]
	\centering
    \includegraphics[width=0.5\linewidth]{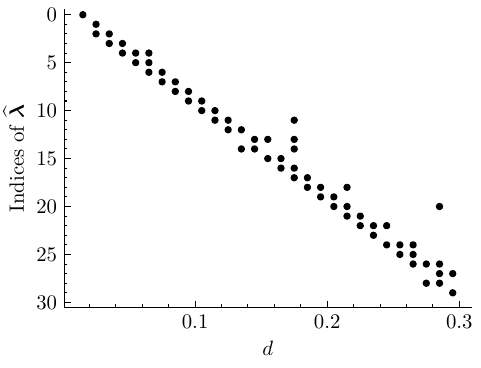}
    \caption{Sparsity of dual reduced dofs $\widehat{\bm{\lambda}}$ selected by greedy active-set method for reconstructions in the validation set of the Hertz problem using a dual dictionary of size 30. Dots indicate the non-zero positions. }
	\label{fig:hertz_dualDof_sparsity}
\end{figure}

\subsection{Application to the Ironing problem}\label{sec:greedy_ironing}
Unlike the Hertz problem, other contact problems might exhibit larger changes in the contact zone, like the ironing problem, where an iron block is pressed against a flat slab and moved along its length (\cref{fig:ironing}). The contact zone varies strongly  with the horizontal displacement parameter $d_x$ (\cref{fig:ironing_snaps}), which worsens the separability~\cite{Kollepara2022}. Therefore, over-complete dictionaries might be more relevant in generating an effective reduced model. 

The ironing problem is considered with two simplifications: the surfaces are assumed to be frictionless and the problem is \emph{quasi-static}, i.e.\ the iron is moved very slowly on the slab surface. The two bodies and surfaces are discretized using same type of elements as the Hertz problem. The iron is discretized using a structured mesh of $15\times15$ nodes, while the slab is discretized using a Cartesian mesh of $10\times50$ nodes.

\begin{figure}[!htb]
	\centering
    \includegraphics[trim=0 25 0 25, clip, width=0.8\linewidth]{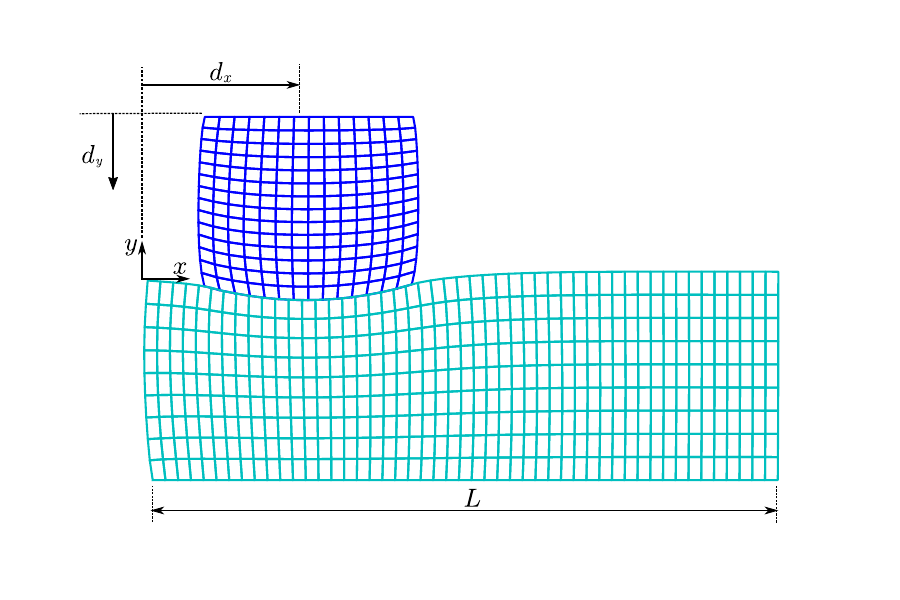}
    \caption{Ironing problem: The iron block is pressed against the flat slab by a displacement $d_y$ and moved horizontally. The horizontal displacement of the iron block $d_x \in [0,L]$ is treated as the parameter in the reduced model}\label{fig:ironing}
\end{figure}

\begin{figure}[!htb]
	\centering
	\includegraphics[width=0.6\linewidth]{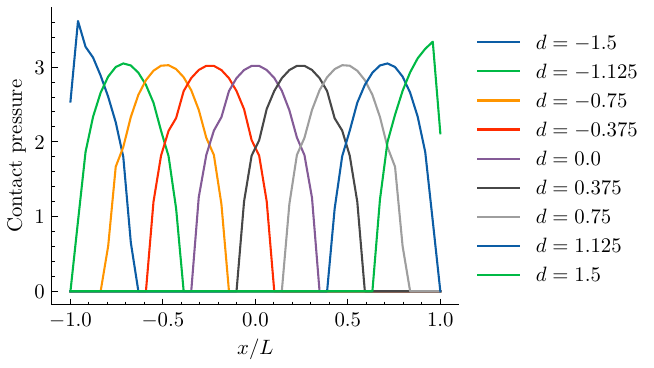}
	\caption{Sample snapshots of contact pressure for ironing problem}\label{fig:ironing_snaps}
\end{figure}

The parametric space is defined as $\mathcal{P}:= \{d_x\ |\ 0\leq d_x\leq L\}$, given length of the slab $L$. The vertical displacement is fixed at $d_y = 0.3$. The training sets $\mathcal{P}_{\texttt{tr}}^n$ are defined in a nested manner, where $n$ is the nested level index. Each training set contains uniformly distributed points in $\mathcal{P}$, such that $\mathcal{P}_{\texttt{tr}}^{n} \subset \mathcal{P}_{\texttt{tr}}^{n+1}$ and $\#\mathcal{P}_{\texttt{tr}}^{n} = 2^{n}+1$. Training sets from nested levels $n=3$ to $n=7$ are used, whereas the validation set is defined as $\mathcal{P}_{\texttt{val}} = \mathcal{P}_{\texttt{tr}}^8 \setminus \mathcal{P}_{\texttt{tr}}^7$. Computational time for each snapshot is around $0.85$ seconds to $1.1$ seconds. 

 The dual dictionary $\DualDict$ is composed of contact pressure snapshots at the $2^n+1$ training points, while the primal basis consists of the left singular vectors truncated by filtering the singular values w.r.t.\ a fixed tolerance $\delta$.

\crefrange{fig:reconstruction_example_ironing}{fig:ironing_hlevel_evolution} contain results of the greedy active-set method applied to the ironing problem. The following are the observations that contrast with the results of Hertz problem:
\begin{enumerate}
    \item \textbf{Inseparability}: \cref{fig:reconstruction_example_ironing} shows the contrast between approximations by dictionaries of different sizes. Its easy to see that the effects of inseparability are more pronounced in the case of ironing problem.
    \item \textbf{Spurious selections}: Comparing \cref{fig:hertz_dualDof_sparsity} and \cref{fig:ironing_dualDof_sparsity}, the sparsity pattern for ironing problem does not exhibit any spurious snapshot selections. This difference can be attributed to the stronger variations in the contact position; unlike the Hertz problem where contact pressure curves were centred at the same position. The activation-deactivation process is, thus, less ambiguous and less sensitive to primal truncations in the case of the ironing problem.
    \item \textbf{Reconstruction errors} in \cref{fig:ironing_hlevel_evolution_errors} are quite similar for $\delta=10^{-10}$ and $\delta=10^{-8}$ possibly because the bottlenecks are the inseparability and finite, although large, size of dual dictionary.
    \item \textbf{Number of iterations and computational time} seem to increase significantly with dictionary size (see \cref{fig:ironing_hlevel_evolution_iters} and \cref{fig:ironing_hlevel_evolution_total_time}). A possible explanation for this behaviour could be that the larger dictionaries of the ironing problem, coupled with stronger parametric dependence of contact zone, force the greedy active-set to perform more iterations to find the appropriate dual dictionary candidates and contact pairs.
\end{enumerate}

\begin{figure}[!htb]
	\centering
	\begin{subfigure}[t]{0.48\linewidth}
        \centering
        \includegraphics[width=1\linewidth]{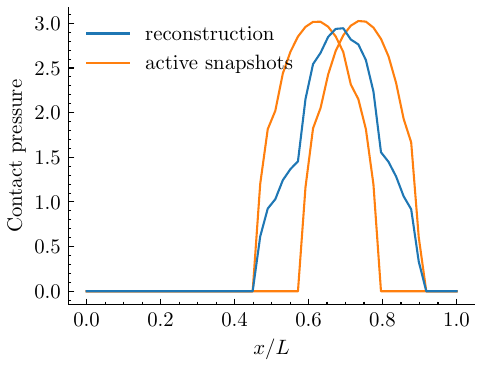}
        \caption{$d_x = 0.68 L$, nested level 3}
	\end{subfigure}
	\begin{subfigure}[t]{0.48\linewidth}
        \centering
        \includegraphics[width=1\linewidth]{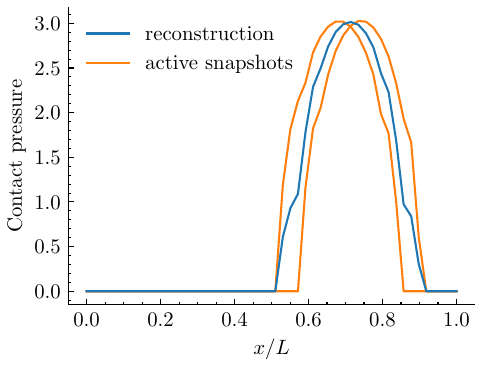}
        \caption{$d_x = 0.65 L$, nested level 4}
	\end{subfigure}
    \caption{Contact pressure reconstructions for certain parametric instances of the ironing problem using greedy active-set method. The dictionary nested level, dictionary elements chosen by the algorithm and the reconstruction errors are also shown.}
	\label{fig:reconstruction_example_ironing}
\end{figure}

\begin{figure}[!htb]
	\centering
    \includegraphics[width=0.5\linewidth]{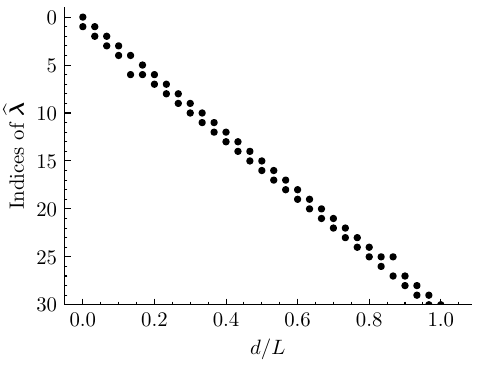}
    \caption{Sparsity of dual reduced dofs $\widehat{\bm{\lambda}}$ selected by greedy active-set method for various reconstructions of ironing problem using the dual dictionary of nested level 4}
	\label{fig:ironing_dualDof_sparsity}
\end{figure}

\begin{figure}[!htpb]
	\centering
    \begin{subfigure}[b]{1.0\linewidth}
        \centering
        \includegraphics[width=0.98\linewidth]{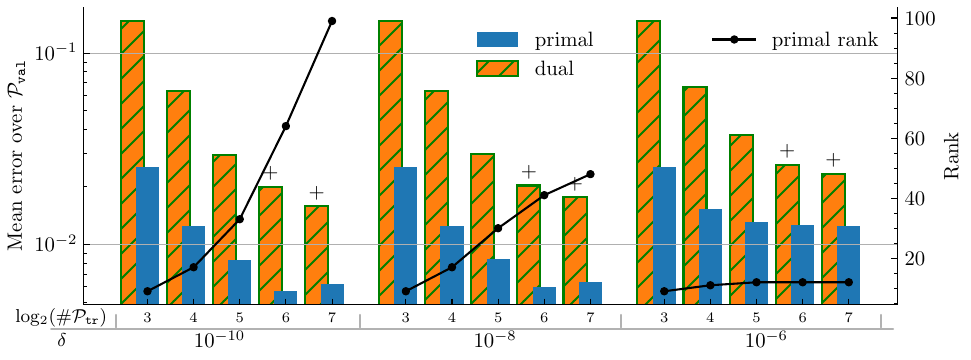}
        \caption{Mean reconstruction errors (bar plot, left $y$-axis) and primal rank (line plot, right $y$-axis). $+$ labels indicate cases where $\# \mathcal{P}_{\text{tr}} > \#$ FOM contact dofs.}
        \label{fig:ironing_hlevel_evolution_errors}
	\end{subfigure}
	\begin{subfigure}[t]{0.48\linewidth}
        \centering
        \includegraphics[width=1\linewidth]{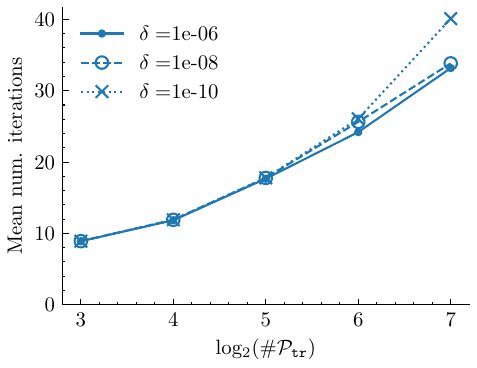}
        \caption{Mean number of iterations}
        \label{fig:ironing_hlevel_evolution_iters}
	\end{subfigure}
	\begin{subfigure}[t]{0.48\linewidth}
        \centering
        \includegraphics[width=1\linewidth]{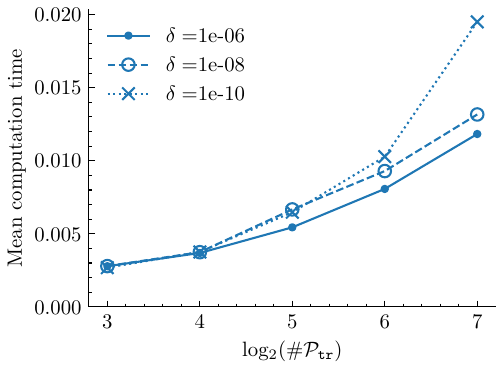}
        \caption{Mean computation time (s), (mean high-fidelity solver time $\approx 0.95$s)}
        \label{fig:ironing_hlevel_evolution_total_time}
	\end{subfigure}
	\begin{subfigure}[b]{0.48\linewidth}
        \centering
        \includegraphics[width=1\linewidth]{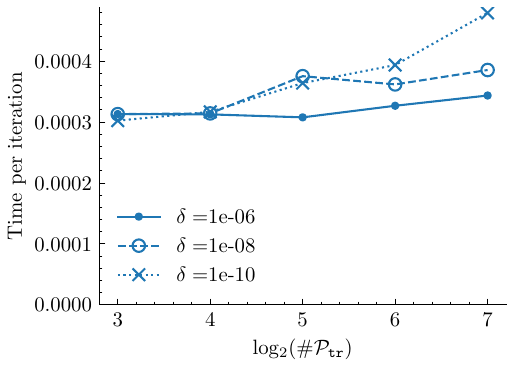}
        \caption{Mean computation time per iteration (s)}
        \label{fig:ironing_hlevel_evolution_time_per_iter}
	\end{subfigure}
    \caption{Evolution of \textbf{(a)} mean reconstruction (relative) errors and primal rank, \textbf{(b)} number of iterations, \textbf{(c)} total and \textbf{(d)} per iteration computation time with training set size for the Ironing problem. Mean values of these quantities over validation set $\mathcal{P}_\texttt{val}$ are shown. Primal and dual errors are computed using $\mathcal{H}^1$ and $\mathcal{L}^2$ norm, respectively. Quantities are shown for different primal truncation tolerance  $\delta \in \{ 10^{6}, \ 10^{-8}$, $10^{-10} \}$. Note that computation time excludes the construction time of non-linear operators. \label{fig:ironing_hlevel_evolution}}	
\end{figure}

\subsection{Application to a two-parameter ironing problem}\label{sec:greedy_ironing_2par}
To validate the same approach on a multidimensional parametric domain, the ironing problem with an additional parameter, i.e.\ the vertical displacement is considered. Thus, the new parametric space is defined as $\mathcal{P} := \{(d_x, d_y)\ |\ 0\leq d_x\leq L, \ 0.1 \leq d_y \leq 0.3 \} $. In the 2D parametric space, the nested training sets are defined such that the number of points is doubled in each parametric direction with increasing levels. The resulting error profiles and computational times are given in \cref{fig:multipar_hlevel_evolution}, with the same trends as the single parameter case.

\begin{figure}[!htpb]
	\centering
    \begin{subfigure}[b]{1.0\linewidth}
        \centering
        \includegraphics[width=0.98\linewidth]{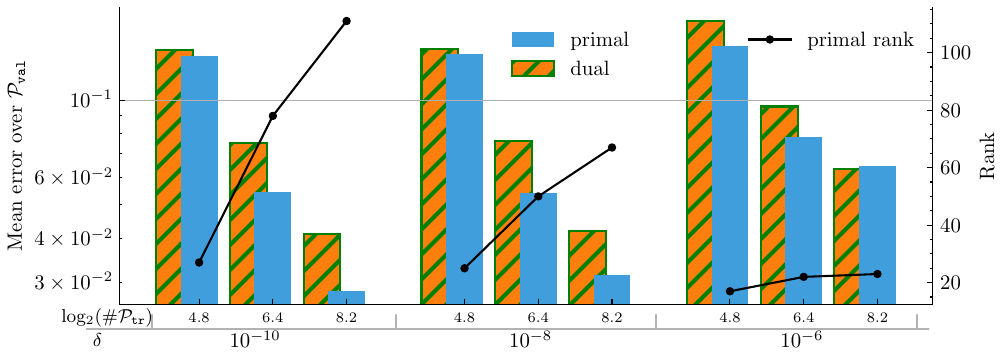}
        \caption{Mean reconstruction errors and primal rank}
	\end{subfigure}
	\begin{subfigure}[b]{0.48\linewidth}
        \centering
        \includegraphics[width=1\linewidth]{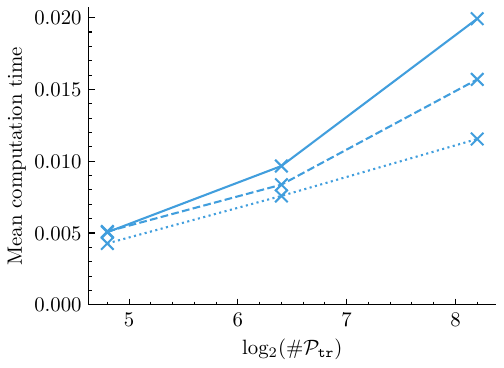}
        \caption{Mean computation time (s)}
        \label{fig:multipar_hlevel_evolution_time}
	\end{subfigure}
	\begin{subfigure}[b]{0.48\linewidth}
        \centering
        \includegraphics[width=1\linewidth]{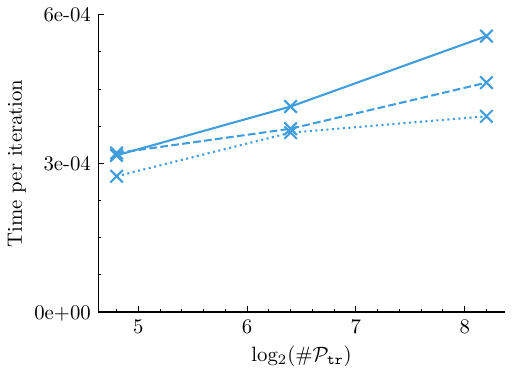}
        \caption{Mean computation time per iteration (s)}
        \label{fig:multipar_hlevel_time_per_iter}
	\end{subfigure}
    \caption{Evolution of \textbf{(a)} mean reconstruction (relative) errors and rank of primal basis, \textbf{(b)} mean computation time and \textbf{(c)} mean computation time per iteration with training set size for the two-parameter ironing problem. \label{fig:multipar_hlevel_evolution}}	
\end{figure}

\section{Discussions}
The application of dictionary methods was proposed to mitigate the challenges posed by linear inseparability of the contact pressure. Unlike a classical reduced basis, an over-complete dictionary has the potential to preserve local information computed in the offline stage, which would otherwise be lost by truncation of a reduced basis. The same has been demonstrated by building reasonably accurate reconstructions using dictionary methods.

The reconstruction errors of the Hertz problem decay consistently when richer training sets are used, and the computational time per iteration seems to be fairly constant, but change slightly with the primal rank. For the ironing problem with large contact variations, the reconstruction errors decays comparatively slower than Hertz problem. Though the dictionary size didn't have any significant influence on the computational time per iteration, the number of iterations (and total time) increased significantly for larger dictionaries in case of ironing problem. This can be explained by the inherently high non-linear behaviour that comes with the large contact surface and could possibly be reduced by using smarter initialization of the contact pairs or by increasing the hierarchical level of the training set gradually.

% Comparing the two dictionary methods: the reconstruction errors of the Hertz problem were better using the greedy active-set algorithm (see \cref{fig:hertz_mean_error_valset,fig:convex_hertz_recon_err}). This observation can be possibly attributed to the use of the monolithic dictionary in case of the convex hull approximations; which forced common interpolation coefficients for the displacement as well as contact pressure. In particular, the convex hull approach shows large errors for small training sets, as the convex hull of a small set of snapshots might be a very small subregion in the entire solution space. On the other hand, the computational time (and number of iterations) for the convex hull approach seems to be less dependent on the dictionary size, compared to the greedy active-set algorithm (see \cref{fig:hertz_num_iters,fig:hertz_computation_time,fig:convex_hertz_computation_time}), thanks to the elimination of non-penetration constraint. In both approaches, the computational time per iteration is nearly flat when the rank is low. Thus, one can infer that large dictionaries and moderate truncation tolerances / rank of $\mathbf{B}$ are an optimal choice. 

% For the convex hull approach, the primary disadvantage is its applicability is limited to a class of contact problems i.e. contact problems with a convex and parameter independent feasible region. The non-convex nature of ironing problem meant that convex hull approximation method could not be used in this case.

As mentioned earlier, the resolution of the inseparabilities and reducing the subsequent errors is not within the scope of this article. The primary objective is to focus on efficient algorithms to utilise over-complete dictionaries for contact problems. A comparison of the greedy active set method with its non-sparse version is shown in Table.~\ref{tab:comparison_methods}. The sparse method provides considerable gains in computational time but not in terms of reconstruction accuracy.

\begin{table}[!htb]
	\setlength\extrarowheight{2pt}
    \small
	\centering
    % \begin{tabular}{lccc}
    \begin{subtable}{.9\linewidth}
	\centering
    \begin{tabular}{@{\extracolsep{6pt}}lcc}
        \toprule
        & \emph{Non-sparse} & \emph{Sparse}  \\
        \cmidrule{2-2} 
        \cmidrule{3-3} 
            % & Active set w.o. CPG & Active set with CPG & Greedy active set  \\
            % & & tolerance $= 1 \times 10^{-2}$ &   \\
            & Active set with CPG & Greedy active set  \\
            & tolerance $= 1 \times 10^{-2}$ &   \\
        \midrule                                                                                              
            % Primal RB size & $22$ & $22$ & $22$ \\
            % Dual RB size & $30$ & $21$ & $30$ \\
            % Time (seconds) & $2.3 \times 10^{-2}$ s & $1.7 \times 10^{-2}$s  & $4 \times 10^{-3}$ s \\
            % Primal errors & $1.7 \times 10^{-3}$ & $1.4 \times 10^{-3}$ & $2 \times 10^{-3}$ \\
            % Dual errors & $5 \times 10^{-2}$ & $5 \times 10^{-2}$  & $5 \times 10^{-2}$ \\
            Primal RB size & $22$ & $22$ \\
            Dual RB/Dictionary size & $21$ & $30$ \\
            Time (seconds) & $1.7 \times 10^{-2}$s  & $4 \times 10^{-3}$ s \\
            Primal reconstruction errors & $1.4 \times 10^{-3}$ & $2 \times 10^{-3}$ \\
            Dual reconstruction errors & $5 \times 10^{-2}$  & $5 \times 10^{-2}$ \\
        \bottomrule
    \end{tabular}
    \caption{Hertz problem with $30$ snapshots in training set.}
	\label{tab:comparison_methods_hertz}
    \vspace{2em}
    \end{subtable}
    \begin{subtable}{.9\linewidth}
	\centering
    \begin{tabular}{@{\extracolsep{6pt}}lcc}
        \toprule
        % & \multicolumn{2}{c}{\emph{Non-sparse}} & \emph{Sparse}  \\
        & \emph{Non-sparse} & \emph{Sparse}  \\
        \cmidrule{2-2} 
        \cmidrule{3-3} 
            % & Active set w.o. CPG & Active set with CPG & Greedy active set  \\
            % & & tolerance $= 5 \times 10^{-2}$ &   \\
            & Active set with CPG & Greedy active set  \\
            & tolerance $= 5 \times 10^{-2}$ &   \\
        \midrule                                                                                              
            % Primal RB size & $99$ & $99$ & $99$ \\
            % Dual RB size & $129$ & $52$ & $129$ \\
            % Time (seconds) & $0.1$ s & $0.08$s  & $0.02$ s \\
            % Primal reconstruction errors & $1.2 \times 10^{-1}$ & $3.7 \times 10^{-2}$ & $3 \times 10^{-2}$ \\
            % Dual reconstruction errors & $8.8 \times 10^{-2}$ & $2.3 \times 10^{-2}$  & $1.5 \times 10^{-2}$ \\
            Primal RB size & $99$ & $99$ \\
            Dual RB/Dictionary size & $52$ & $129$ \\
            Time (seconds) & $0.08$s  & $0.02$ s \\
            Primal reconstruction errors & $3.7 \times 10^{-2}$ & $3 \times 10^{-2}$ \\
            Dual reconstruction errors & $2.3 \times 10^{-2}$  & $1.5 \times 10^{-2}$ \\
        \bottomrule
    \end{tabular}
    \caption{Ironing problem with $129$ snapshots in training set.}
	\label{tab:comparison_methods_ironing}
    \end{subtable}
    \caption{Performance comparison of non-sparse and sparse methods. Mean of the quantities (time and reconstruction errors) observed on the validation set are given.}
	\label{tab:comparison_methods}
\end{table}

\section{Conclusions and Perspectives}
\emph{Conclusions}: The key idea of using over-complete dictionaries for contact mechanics problems is introduced, whose necessity is justified using the known limitations of current reduced methods. Based on the concepts of sparse regressions, an algorithm is proposed to solve dictionary-based approximations of contact problems. Approximations based on large dictionaries were shown to be limit the inseparability issues to a certain extent. Performance of the proposed method seems to be correlated to the contact position variability, with satisfactory performance for small variabilities.

\emph{Perspectives}: The primary bottleneck in the proposed method is the evaluation of non-linear operators that describe the inequality constraint i.e.\ the non-penetration condition. Current methods also run into issues caused by local properties of contact, as contact pairs change significantly for problems with large contact variations. Moreover, non-linear methods to interpolate the solution manifold should be explored.

\section*{Acknowledgements}
The research has been carried out as part of a doctoral thesis funded by IRT Jules Verne under the PERFORM programme. All computations were done on \texttt{python} using open source packages, namely NumPy~\cite{numpy} and SciPy~\cite{scipy} for scientific computations and Matplotlib~\cite{matplotlib} for creating plots. Meshes were created using Gmsh~\cite{gmsh}.
%%fakesection bibliography
\bibliography{bibfile}

\begin{thebibliography}{52}
\providecommand{\natexlab}[1]{#1}
\providecommand{\url}[1]{\texttt{#1}}
\expandafter\ifx\csname urlstyle\endcsname\relax
  \providecommand{\doi}[1]{doi: #1}\else
  \providecommand{\doi}{doi: \begingroup \urlstyle{rm}\Url}\fi

\bibitem[Li et~al.(2022)Li, Li, Zhou, Feng, and Zhou]{Li2022}
B.~Li, P.~Li, R.~Zhou, X.-Q. Feng, and K.~Zhou.
\newblock Contact mechanics in tribological and contact damage-related
  problems: A review.
\newblock \emph{Tribology International}, 171:\penalty0 107534, 2022.
\newblock ISSN 0301-679X.
\newblock \doi{10.1016/j.triboint.2022.107534}.

\bibitem[Brown et~al.(2000)Brown, Attaway, Plimpton, and
  Hendrickson]{Brown2000}
K.~Brown, S.~Attaway, S.~Plimpton, and B.~Hendrickson.
\newblock Parallel strategies for crash and impact simulations.
\newblock \emph{Computer Methods in Applied Mechanics and Engineering},
  184\penalty0 (2):\penalty0 375--390, 2000.
\newblock ISSN 0045-7825.
\newblock \doi{10.1016/S0045-7825(99)00235-2}.

\bibitem[Button(2013)]{Button2013}
S.~T. Button.
\newblock \emph{Tribology in Metal Forming Processes}, pages 103--120.
\newblock Springer Berlin Heidelberg, Berlin, Heidelberg, 2013.
\newblock ISBN 978-3-642-31683-8.
\newblock \doi{10.1007/978-3-642-31683-8_3}.

\bibitem[Ateshian et~al.(2015)Ateshian, Henak, and Weiss]{Ateshian2015}
G.~A. Ateshian, C.~R. Henak, and J.~A. Weiss.
\newblock Toward patient-specific articular contact mechanics.
\newblock \emph{Journal of Biomechanics}, 48\penalty0 (5):\penalty0 779--786,
  2015.
\newblock ISSN 0021-9290.
\newblock \doi{10.1016/j.jbiomech.2014.12.020}.

\bibitem[Wriggers and Simo(1985)]{Wriggers1985}
P.~Wriggers and J.~C. Simo.
\newblock A note on tangent stiffness for fully nonlinear contact problems.
\newblock \emph{Communications in Applied Numerical Methods}, 1\penalty0
  (5):\penalty0 199--203, 1985.
\newblock \doi{10.1002/cnm.1630010503}.

\bibitem[Hallquist et~al.(1985)Hallquist, Goudreau, and Benson]{Hallquist1985}
J.~Hallquist, G.~Goudreau, and D.~Benson.
\newblock Sliding interfaces with contact-impact in large-scale lagrangian
  computations.
\newblock \emph{Computer Methods in Applied Mechanics and Engineering},
  51\penalty0 (1):\penalty0 107--137, 1985.
\newblock \doi{10.1016/0045-7825(85)90030-1}.

\bibitem[Fischer and Wriggers(2006)]{Fischer2006}
K.~A. Fischer and P.~Wriggers.
\newblock {Mortar based frictional contact formulation for higher order
  interpolations using the moving friction cone}.
\newblock \emph{Computer Methods in Applied Mechanics and Engineering},
  195\penalty0 (37-40):\penalty0 5020--5036, 2006.
\newblock \doi{10.1016/j.cma.2005.09.025}.

\bibitem[Wriggers(2006)]{Wriggers2006}
P.~Wriggers.
\newblock \emph{{Computational Contact Mechanics}}.
\newblock Springer Berlin Heidelberg, 2006.
\newblock ISBN 978-3-540-32608-3.
\newblock \doi{10.1007/978-3-540-32609-0}.

\bibitem[Chouly et~al.(2017)Chouly, Fabre, Hild, Mlika, Pousin, and
  Renard]{Chouly2017}
F.~Chouly, M.~Fabre, P.~Hild, R.~Mlika, J.~Pousin, and Y.~Renard.
\newblock {An Overview of Recent Results on Nitsche's Method for Contact
  Problems}.
\newblock pages 93--141. 2017.
\newblock \doi{10.1007/978-3-319-71431-8_4}.

\bibitem[Yastrebov(2013)]{Yastrebov2013}
V.~A. Yastrebov.
\newblock \emph{{Numerical Methods in Contact Mechanics}}.
\newblock John Wiley \& Sons, Inc., Hoboken, NJ USA, feb 2013.
\newblock ISBN 9781118647974.
\newblock \doi{10.1002/9781118647974}.

\bibitem[Meier et~al.(2005)Meier, López, Monserrat, Juan, and
  Alcañiz]{Meier2005}
U.~Meier, O.~López, C.~Monserrat, M.~Juan, and M.~Alcañiz.
\newblock Real-time deformable models for surgery simulation: a survey.
\newblock \emph{Computer Methods and Programs in Biomedicine}, 77\penalty0
  (3):\penalty0 183--197, 2005.
\newblock ISSN 0169-2607.
\newblock \doi{https://doi.org/10.1016/j.cmpb.2004.11.002}.
\newblock URL
  \url{https://www.sciencedirect.com/science/article/pii/S0169260704002093}.

\bibitem[Tan et~al.(2023)Tan, Chen, Zou, Yang, and Du]{Tan2023}
X.~Tan, W.~Chen, T.~Zou, J.~Yang, and B.~Du.
\newblock Real-time prediction of mechanical behaviors of underwater shield
  tunnel structure using machine learning method based on structural health
  monitoring data.
\newblock \emph{Journal of Rock Mechanics and Geotechnical Engineering},
  15\penalty0 (4):\penalty0 886--895, 2023.
\newblock ISSN 1674-7755.
\newblock \doi{https://doi.org/10.1016/j.jrmge.2022.06.015}.

\bibitem[Haasdonk et~al.(2012)Haasdonk, Salomon, and Wohlmuth]{Haasdonk2012}
B.~Haasdonk, J.~Salomon, and B.~Wohlmuth.
\newblock {A Reduced Basis Method for Parametrized Variational Inequalities}.
\newblock \emph{SIAM J. Numer. Anal.}, 50\penalty0 (5):\penalty0 2656--2676,
  jan 2012.
\newblock \doi{10.1137/110835372}.

\bibitem[Balajewicz et~al.(2016)Balajewicz, Amsallem, and
  Farhat]{Balajewicz2016}
M.~Balajewicz, D.~Amsallem, and C.~Farhat.
\newblock {Projection-based model reduction for contact problems}.
\newblock \emph{Int. J. Numer. Methods Eng.}, 106\penalty0 (8):\penalty0
  644--663, may 2016.
\newblock \doi{10.1002/nme.5135}.

\bibitem[Niakh et~al.(2023)Niakh, Drouet, Ehrlacher, and Ern]{Niakh_mCPG_2023}
I.~Niakh, G.~Drouet, V.~Ehrlacher, and A.~Ern.
\newblock Stable model reduction for linear variational inequalities with
  parameter-dependent constraints.
\newblock \emph{ESAIM: Mathematical Modelling and Numerical Analysis},
  57\penalty0 (1):\penalty0 167--189, Jan. 2023.
\newblock ISSN 2822-7840, 2804-7214.
\newblock \doi{10.1051/m2an/2022077}.

\bibitem[Benaceur et~al.(2020)Benaceur, Ern, and Ehrlacher]{Benaceur2020}
A.~Benaceur, A.~Ern, and V.~Ehrlacher.
\newblock {A reduced basis method for parametrized variational inequalities
  applied to contact mechanics}.
\newblock \emph{Int. J. Numer. Methods Eng.}, 121\penalty0 (6):\penalty0
  1170--1197, mar 2020.
\newblock \doi{10.1002/nme.6261}.

\bibitem[Fauque et~al.(2018)Fauque, Rami{\`{e}}re, and Ryckelynck]{Fauque2018}
J.~Fauque, I.~Rami{\`{e}}re, and D.~Ryckelynck.
\newblock {Hybrid hyper-reduced modeling for contact mechanics problems}.
\newblock \emph{Int. J. Numer. Methods Eng.}, 115\penalty0 (1):\penalty0
  117--139, Jul 2018.
\newblock \doi{10.1002/nme.5798}.

\bibitem[Le~Berre et~al.(2022)Le~Berre, Ramière, Fauque, and
  Ryckelynck]{LeBerre2022}
S.~Le~Berre, I.~Ramière, J.~Fauque, and D.~Ryckelynck.
\newblock Condition {Number} and {Clustering}-{Based} {Efficiency}
  {Improvement} of {Reduced}-{Order} {Solvers} for {Contact} {Problems} {Using}
  {Lagrange} {Multipliers}.
\newblock \emph{Mathematics}, 10\penalty0 (9):\penalty0 1495, Apr. 2022.
\newblock ISSN 2227-7390.
\newblock \doi{10.3390/math10091495}.

\bibitem[Bader et~al.(2016)Bader, Zhang, and Veroy]{Bader2016}
E.~Bader, Z.~Zhang, and K.~Veroy.
\newblock {An empirical interpolation approach to reduced basis approximations
  for variational inequalities}.
\newblock \emph{Math. Comput. Model. Dyn. Syst.}, 22\penalty0 (4):\penalty0
  345--361, 2016.
\newblock \doi{10.1080/13873954.2016.1198388}.

\bibitem[Niakh et~al.(2023)Niakh, Drouet, Ehrlacher, and
  Ern]{Niakh_Nitsche_2023}
I.~Niakh, G.~Drouet, V.~Ehrlacher, and A.~Ern.
\newblock A reduced basis method for frictional contact problems formulated
  with {Nitsche}'s method.
\newblock 2023.
\newblock \doi{10.48550/ARXIV.2307.11541}.
\newblock Publisher: arXiv Version Number: 1.

\bibitem[Manvelyan et~al.(2021)Manvelyan, Simeon, and Wever]{Manvelyan2021}
D.~Manvelyan, B.~Simeon, and U.~Wever.
\newblock {An efficient model order reduction scheme for dynamic contact in
  linear elasticity}.
\newblock \emph{Computational Mechanics}, 68\penalty0 (6):\penalty0 1283--1295,
  dec 2021.
\newblock \doi{10.1007/s00466-021-02068-4}.

\bibitem[Mulye et~al.(2020)Mulye, Hemmer, Morançay, Binetruy, Leygue,
  Comas-Cardona, Pichon, and Guillon]{Mulye2020}
P.~Mulye, J.~Hemmer, L.~Morançay, C.~Binetruy, A.~Leygue, S.~Comas-Cardona,
  P.~Pichon, and D.~Guillon.
\newblock Numerical modeling of interply adhesion in composite forming of
  viscous discontinuous thermoplastic prepregs.
\newblock \emph{Composites Part B: Engineering}, 191:\penalty0 107953, 2020.
\newblock ISSN 1359-8368.
\newblock \doi{10.1016/j.compositesb.2020.107953}.

\bibitem[Kollepara et~al.(2022)Kollepara, Navarro‐Jiménez, Le~Guennec,
  Silva, and Aguado]{Kollepara2022}
K.~S. Kollepara, J.~M. Navarro‐Jiménez, Y.~Le~Guennec, L.~Silva, and J.~V.
  Aguado.
\newblock On the limitations of low‐rank approximations in contact mechanics
  problems.
\newblock \emph{International Journal for Numerical Methods in Engineering},
  page nme.7118, Sept. 2022.
\newblock ISSN 0029-5981, 1097-0207.
\newblock \doi{10.1002/nme.7118}.

\bibitem[Davies(2002)]{davies2002}
R.~H. Davies.
\newblock \emph{{Learning shape: optimal models for analysing natural
  variability}}.
\newblock PhD thesis, University of Manchester, 2002.
\newblock URL \url{http://www.bmva.org/theses/2002/2002-davies.pdf}.

\bibitem[Amsallem et~al.(2012)Amsallem, Zahr, and Farhat]{Amsallem2012}
D.~Amsallem, M.~J. Zahr, and C.~Farhat.
\newblock Nonlinear model order reduction based on local reduced-order bases.
\newblock \emph{International Journal for Numerical Methods in Engineering},
  92\penalty0 (10):\penalty0 891--916, 2012.
\newblock \doi{10.1002/nme.4371}.

\bibitem[Amsallem et~al.(2009)Amsallem, Cortial, Carlberg, and
  Farhat]{Amsallem2009}
D.~Amsallem, J.~Cortial, K.~Carlberg, and C.~Farhat.
\newblock A method for interpolating on manifolds structural dynamics
  reduced-order models.
\newblock \emph{International Journal for Numerical Methods in Engineering},
  80\penalty0 (9):\penalty0 1241--1258, 2009.
\newblock \doi{10.1002/nme.2681}.

\bibitem[Aharon et~al.(2006)Aharon, Elad, and Bruckstein]{Aharon2006}
M.~Aharon, M.~Elad, and A.~Bruckstein.
\newblock K-svd: An algorithm for designing overcomplete dictionaries for
  sparse representation.
\newblock \emph{IEEE Transactions on Signal Processing}, 54\penalty0
  (11):\penalty0 4311–4322, 2006.
\newblock \doi{10.1109/tsp.2006.881199}.

\bibitem[Rubinstein et~al.(2014)Rubinstein, Zibulevsky, and
  Elad]{Rubinstein2014}
R.~Rubinstein, M.~Zibulevsky, and M.~Elad.
\newblock {Efficient Implementation of the K-SVD Algorithm Using Batch
  Orthogonal Matching Pursuit}.
\newblock \emph{CS Techion}, \penalty0 (January 2008), 2014.

\bibitem[Balabanov(2021)]{Balabanov2021}
A.~Balabanov, O.and~Nouy.
\newblock {Randomized linear algebra for model reduction—part II: minimal
  residual methods and dictionary-based approximation}.
\newblock \emph{Advances in Computational Mathematics}, 47\penalty0
  (2):\penalty0 26, apr 2021.
\newblock \doi{10.1007/s10444-020-09836-5}.

\bibitem[Lauzeral(2019)]{Lauzeral2019}
N.~Lauzeral.
\newblock \emph{Reduced order and sparse representations for patient-specific
  modeling in computational surgery}.
\newblock PhD thesis, {{\'E}cole centrale de Nantes}, dec 2019.
\newblock URL \url{http://www.theses.fr/2019ECDN0062}.

\bibitem[Simoncelli et~al.(1992)Simoncelli, Freeman, Adelson, and
  Heeger]{Simoncelli1992}
E.~Simoncelli, W.~Freeman, E.~Adelson, and D.~Heeger.
\newblock Shiftable multiscale transforms.
\newblock \emph{IEEE Transactions on Information Theory}, 38\penalty0
  (2):\penalty0 587--607, 1992.
\newblock \doi{10.1109/18.119725}.

\bibitem[Cand{\`e}s and Donoho(2002)]{Candes2002}
E.~J. Cand{\`e}s and D.~L. Donoho.
\newblock {Recovering edges in ill-posed inverse problems: optimality of
  curvelet frames}.
\newblock \emph{The Annals of Statistics}, 30\penalty0 (3):\penalty0 784 --
  842, 2002.
\newblock \doi{10.1214/aos/1028674842}.
\newblock URL \url{https://doi.org/10.1214/aos/1028674842}.

\bibitem[Elad(2010)]{Elad2010}
M.~Elad.
\newblock \emph{{Sparse and Redundant Representations}}.
\newblock Springer New York, New York, NY, 2010.
\newblock ISBN 978-1-4419-7010-7.
\newblock \doi{10.1007/978-1-4419-7011-4}.

\bibitem[Brunton et~al.(2016)Brunton, Proctor, and Kutz]{Brunton2016}
S.~L. Brunton, J.~L. Proctor, and J.~N. Kutz.
\newblock {Discovering governing equations from data by sparse identification
  of nonlinear dynamical systems}.
\newblock \emph{Proceedings of the National Academy of Sciences}, 113\penalty0
  (15):\penalty0 3932--3937, apr 2016.
\newblock \doi{10.1073/pnas.1517384113}.

\bibitem[Mairal et~al.(2014)Mairal, Bach, and Ponce]{Mairal2014}
J.~Mairal, F.~Bach, and J.~Ponce.
\newblock {Sparse modeling for image and vision processing}.
\newblock \emph{Foundations and Trends in Computer Graphics and Vision},
  8\penalty0 (2-3):\penalty0 85--283, 2014.
\newblock \doi{10.1561/0600000058}.

\bibitem[Tibshirani(1996)]{Tibshirani1996}
R.~Tibshirani.
\newblock {Regression Shrinkage and Selection Via the Lasso}.
\newblock \emph{Journal of the Royal Statistical Society: Series B
  (Methodological)}, 58\penalty0 (1):\penalty0 267--288, jan 1996.
\newblock \doi{10.1111/j.2517-6161.1996.tb02080.x}.

\bibitem[Candes and Tao(2007)]{Candes2007}
E.~Candes and T.~Tao.
\newblock {The Dantzig selector: Statistical estimation when p is much larger
  than n}.
\newblock \emph{The Annals of Statistics}, 35\penalty0 (6):\penalty0
  2313--2351, dec 2007.
\newblock \doi{10.1214/009053606000001523}.

\bibitem[Pati et~al.(1993)Pati, Rezaiifar, and Krishnaprasad]{Pati1993}
Y.~Pati, R.~Rezaiifar, and P.~Krishnaprasad.
\newblock Orthogonal matching pursuit: recursive function approximation with
  applications to wavelet decomposition.
\newblock In \emph{Proceedings of 27th Asilomar Conference on Signals, Systems
  and Computers}, pages 40--44 vol.1, 1993.
\newblock \doi{10.1109/ACSSC.1993.342465}.

\bibitem[Gorodnitsky and Rao(1997)]{Gorodnitsky1997}
I.~Gorodnitsky and B.~Rao.
\newblock {Sparse signal reconstruction from limited data using FOCUSS: a
  re-weighted minimum norm algorithm}.
\newblock \emph{IEEE Transactions on Signal Processing}, 45\penalty0
  (3):\penalty0 600--616, mar 1997.
\newblock \doi{10.1109/78.558475}.

\bibitem[Wang et~al.(2011)Wang, Yang, Lai, Kovanis, and Grebogi]{Wang2011}
W.-X. Wang, R.~Yang, Y.-C. Lai, V.~Kovanis, and C.~Grebogi.
\newblock Predicting catastrophes in nonlinear dynamical systems by compressive
  sensing.
\newblock \emph{Phys. Rev. Lett.}, 106:\penalty0 154101, Apr 2011.
\newblock \doi{10.1103/PhysRevLett.106.154101}.

\bibitem[Schaeffer et~al.(2013)Schaeffer, Caflisch, Hauck, and
  Osher]{Schaeffer2013}
H.~Schaeffer, R.~Caflisch, C.~D. Hauck, and S.~Osher.
\newblock Sparse dynamics for partial differential equations.
\newblock \emph{Proceedings of the National Academy of Sciences}, 110\penalty0
  (17):\penalty0 6634--6639, 2013.
\newblock \doi{10.1073/pnas.1302752110}.

\bibitem[Halko et~al.(2011)Halko, Martinsson, and Tropp]{Halko2011}
N.~Halko, P.~G. Martinsson, and J.~A. Tropp.
\newblock Finding structure with randomness: Probabilistic algorithms for
  constructing approximate matrix decompositions.
\newblock \emph{SIAM Review}, 53\penalty0 (2):\penalty0 217--288, 2011.
\newblock \doi{10.1137/090771806}.

\bibitem[Le~Guennec et~al.(2018)Le~Guennec, Brunet, Daim, Chau, and
  Tourbier]{LeGuennec2018}
Y.~Le~Guennec, J.~Brunet, F.~Daim, M.~Chau, and Y.~Tourbier.
\newblock A parametric and non-intrusive reduced order model of car crash
  simulation.
\newblock \emph{Computer Methods in Applied Mechanics and Engineering},
  338:\penalty0 186 -- 207, 2018.
\newblock \doi{10.1016/j.cma.2018.03.005}.

\bibitem[Harris et~al.(2020)Harris, Millman, van~der Walt, and et. al]{numpy}
C.~R. Harris, K.~J. Millman, S.~J. van~der Walt, and et. al.
\newblock Array programming with {NumPy}.
\newblock \emph{Nature}, 585\penalty0 (7825):\penalty0 357--362, Sept. 2020.
\newblock \doi{10.1038/s41586-020-2649-2}.

\bibitem[Virtanen et~al.(2020)Virtanen, Gommers, Oliphant, and et. al]{scipy}
P.~Virtanen, R.~Gommers, T.~Oliphant, and et. al.
\newblock {{SciPy} 1.0: Fundamental Algorithms for Scientific Computing in
  Python}.
\newblock \emph{Nature Methods}, 17:\penalty0 261--272, 2020.
\newblock \doi{10.1038/s41592-019-0686-2}.

\bibitem[Hunter(2007)]{matplotlib}
J.~D. Hunter.
\newblock Matplotlib: A 2d graphics environment.
\newblock \emph{Computing in Science \& Engineering}, 9\penalty0 (3):\penalty0
  90--95, 2007.
\newblock \doi{10.1109/MCSE.2007.55}.

\bibitem[Geuzaine and Remacle(2008)]{gmsh}
C.~Geuzaine and J.-F. Remacle.
\newblock Gmsh: A 3-d finite element mesh generator with built-in pre- and
  post-processing facilities.
\newblock \emph{International Journal for Numerical Methods in Engineering},
  79:\penalty0 1309--1331, 2008.

\bibitem[Mallat and Zhang(1993)]{Mallat1992}
S.~Mallat and Z.~Zhang.
\newblock Matching pursuits with time-frequency dictionaries.
\newblock \emph{IEEE Transactions on Signal Processing}, 41\penalty0
  (12):\penalty0 3397--3415, 1993.
\newblock \doi{10.1109/78.258082}.

\bibitem[Manucci et~al.(2022)Manucci, Aguado, and Borzacchiello]{Manucci2022}
M.~Manucci, J.~V. Aguado, and D.~Borzacchiello.
\newblock Sparse data-driven quadrature rules via $\ell^p$-quasi-norm
  minimization.
\newblock \emph{Computational Methods in Applied Mathematics}, 2022.
\newblock \doi{10.1515/cmam-2021-0131}.

\bibitem[Lawson and Hanson(1995)]{Lawson1995}
C.~L. Lawson and R.~J. Hanson.
\newblock \emph{{Solving Least Squares Problems}}.
\newblock Society for Industrial and Applied Mathematics, jan 1995.
\newblock ISBN 978-0-89871-356-5.
\newblock \doi{10.1137/1.9781611971217}.

\bibitem[Veroy and Patera(2005)]{Veroy2005}
K.~Veroy and A.~T. Patera.
\newblock Certified real-time solution of the parametrized steady
  incompressible navier–stokes equations: rigorous reduced-basis a posteriori
  error bounds.
\newblock \emph{International Journal for Numerical Methods in Fluids},
  47\penalty0 (8-9):\penalty0 773--788, 2005.
\newblock \doi{10.1002/fld.867}.

\bibitem[Liberge and Hamdouni(2010)]{Liberge2010}
E.~Liberge and A.~Hamdouni.
\newblock Reduced order modelling method via proper orthogonal decomposition
  (pod) for flow around an oscillating cylinder.
\newblock \emph{Journal of Fluids and Structures}, 26\penalty0 (2):\penalty0
  292--311, 2010.
\newblock ISSN 0889-9746.
\newblock \doi{10.1016/j.jfluidstructs.2009.10.006}.

\end{thebibliography}

\newpage 
\appendix
\section{Sparse methods}\label{sec:sparse_methods_list}
As discussed in \cref{sec:sparse_methods}, solving the problem~\eqref{eq:parsimony_idea} with $\ell_0$-norm is NP-hard. Thus, practically usable sparse regression methods have been devised by the community, such as, sparsity inducing surrogate norms in place of the $\ell_0$-norm, or greedy methods. Some examples include LASSO~\cite{Tibshirani1996} and {Dantzig Selector}~\cite{Candes2007} that uses $\ell_0$-norm. The methods used in this work, OMP and FOCUSS are detailed below:

\subsection{Orthogonal Matching Pursuit}
The OMP method~\cite{Pati1993} aims to approximate the solution of~\eqref{eq:parsimony_idea} with $p=0$ using a greedy technique. At each greedy step $k$, the OMP algorithm searches for the column of dictionary $\mathbf{D}$ that is the least orthogonal to the current residual $(\bm{x}-\mathbf{D} \bm{\alpha}^{k-1})$ and adds to the set of indices previously selected. This process is carried out until the residual is below the tolerance $\varepsilon$. This is done by projecting the residual vector on each column and selecting the largest projection.

\begin{align*}
    \mathcal{I} \gets \mathcal{I} \cup \argmax_j \left |\mathbf{D}[:,j]^T (\bm{x}-\mathbf{D} \bm{\alpha}^{k-1})\right |
    % \label{eq:omp_identification}
\end{align*}

where $\mathcal{I}$ is the current set of selected indices and columns of dictionary $\mathbf{D}$ are normalized. Then the new set of coefficients $\bm{\alpha}^k$ are computed using the least square solution:

\begin{align*}
        \bm{\alpha}^k_\mathcal{I} = (\mathbf{D}_{\mathcal{I}}^T \mathbf{D}_{\mathcal{I}})^{-1} \mathbf{D}_{\mathcal{I}}^T\bm{x}
        \label{eq:omp_updaterep}
\end{align*}

with $\bm{\alpha}_\mathcal{I} =\bm{\alpha}[\mathcal{I}]$ and $\mathbf{D}_\mathcal{I} = \mathbf{D}[:,\mathcal{I}]$.

Update of $\mathcal{I}$ is locally optimal, but the above least-square update of $\bm{\alpha}^k$, the OMP solution is optimal w.r.t to the currently selected subset $\mathbf{D}_\mathcal{I}$ of the dictionary. This is in contrast to the predecessor of OMP, the Matching Pursuit~\cite{Mallat1992}, which used locally optimal updates for $\bm{\alpha}^k$ at the current iteration.

\subsection{FOCUSS}
The FOCUSS method~\cite{Gorodnitsky1997} finds a sparse approximation using an iterative process that starts from a fully dense solution and progresses towards ``localized energy solutions''. The FOCUSS method is broadly based on the following ideas:
\begin{itemize}
    \item For a full rank dictionary of size, say $m \times n\ (m < n)$, the system $\mathbf{D} \bm{\alpha} = \bm{x}$ is under-determined. For such a system, the closed form solution for~\eqref{eq:parsimony_idea} with $p=2$ is given by $\bm{\alpha} = \mathbf{D}^\dagger \bm{x}$, where $\mathbf{D}^\dagger$ is the Moore-Penrose pseudo-inverse. Note that this solution is not sparse in general, as $\ell_2$-norm does not have sparsity-inducing properties. 
    \item A modified minimization problem can be defined by replacing the $\norm{\bm{\alpha}}{2}$ with a weighted norm $\norm{\mathbf{W}^{-1} \bm{\alpha}}{2}$, or more generically $\norm{\mathbf{W}^\dagger \bm{\alpha}}{2}$ if $\mathbf{W}$ is singular. The authors of~\cite{Gorodnitsky1997} argue that by changing $\mathbf{W}$, every possible solution of the under-determined system $\mathbf{D} \bm{\alpha} = \bm{x}$ can be obtained.
    \item The basis of the algorithm lies in iterating towards a sparse solution using a weight $\mathbf{W}^k$ which induces sparsity. For the iteration $k$, the weight is chosen as $\mathbf{W}_{k} = \mathtt{diag}(\bm{\alpha}^{k-1})$. The trick lies in the fact that the algorithm, effectively, minimizes the following weighted norm:
        \begin{equation*}
            \norm{\mathbf{W}_k^\dagger \bm{\alpha}}{2}^2 = \sum_i \left (\frac{\alpha^k_i}{\alpha^{k-1}_i} \right )^2 
        \end{equation*}
    due to which smaller entries of $\alpha^{k-1}$ tend to diminish further.
\end{itemize}
On combining these ideas, the FOCUSS algorithm ends up with quite a simple implementation, given in \cref{alg:focuss}. Usually, the algorithm is initialized with the $\bm{\alpha}^0$ containing all non-zeros. In~\cite{Gorodnitsky1997}, the $\norm{\bm{\alpha}}{2}$ minimizing solution $\bm{\alpha}^0 = \mathbf{D}^\dagger \bm{x}$ was used for initialization. 
\begin{algorithm}[h]
    \caption{FOCUSS} \label{alg:focuss}
    \begin{algorithmic}[1]
        \State{Inputs: $\mathbf{D}, \bm{x}$, $\bm{\alpha}^0$}
        \State{Initialize $k = 1$}
        % \Statex{\phantom{Initialize }$\bm{\alpha}^0 = \mathbf{D}^\dagger \bm{x} $}
        \While{$\bm{\alpha}$ not converge}
        \State{$\mathbf{W}^k = \mathtt{diag}(\bm{\alpha}^{k-1})$}
        \State{$\bm{\alpha}^k = \mathbf{W}^k (\mathbf{D} \mathbf{W}^k)^\dagger \bm{x}$} \label{step:focuss_alpha}
        \State{$k \gets k+1$}
        \EndWhile{}
    \end{algorithmic}
\end{algorithm}

\noindent The \textbf{nnFOCUSS} Algorithm: The FOCUSS algorithm can be modified to generate a non-negative solution~\cite{Lauzeral2019,Manucci2022}. The nnFOCUSS algorithm works by computing an appropriate relaxation parameter for each iterative update that maintains non-negativity. After computation of new coefficients at Step~\ref{step:focuss_alpha} of the \cref{alg:focuss}, the relaxation step of nnFOCUSS algorithm is performed as follows: \\

\begin{algorithmic}[0]
    \If{$\min(\bm{\alpha}^k) < 0$}
    \State{$\Delta \bm{\alpha} = -(\bm{\alpha}^k - \bm{\alpha}^{k-1})^{-}$}
    \State{$\bm{\alpha}^k \gets \bm{\alpha}^{k-1} + \min \left (\frac{\bm{\alpha}^{k-1}}{\Delta \bm{\alpha}} \right ) \Delta \bm{\alpha}$} \Comment{Element wise division}
    \EndIf{}
\end{algorithmic}

\noindent The algorithm must be initialized using non-negative coefficients, $\bm{\alpha}^0$, which is computed using the non-negative least squares~\cite{Lawson1995} solution (which is not sparse, in general). 
\section{A non-penetrating convex hull approach with monolithic dictionaries}\label{sec:convex_hull}
% Two dictionary-based approximation methods for contact mechanics problems have been formulated, namely the \emph{convex-hull approximation} approach and \emph{greedy active-set} approach. To provide a broad overview of the features of these methods, \cref{tab:dictionary_methods} gives a list of features of these methods. The important differences between the two methods are bases types and the manner of enforcing constraints. The two formulations and their application to contact problems are discussed in this section.
Another novel sparse approximation method was formulated, namely the \emph{convex-hull approximation} approach. This section briefly addresses the formulation and preliminary results. The convex hull method differs from the greedy active-set method in several aspects. Before this method is formulated, a broad overview of these differences is provided in \cref{tab:dictionary_methods}.

\begin{table}[!htb]
	\setlength\extrarowheight{2pt}
    \small
	\centering
    % \begin{tabular}{lcc}
    \begin{tabular}{m{0.2\textwidth}m{0.35\textwidth}m{0.05\textwidth}m{0.31\textwidth}}
        \toprule
                                & Convex hull approximation                         & & Greedy active-set       \\
        \midrule                                                                                              
        Primal Basis type       & Dictionary                                        & & Truncated low-rank      \\
        Dual Basis type         & Dictionary                                        & & Dictionary              \\
        Monolithic dictionary   & Yes                                               & & No                      \\
        Iteration method        & Fixed point                                       & & Fixed point             \\ 
        Constraint enforcement  & Convex combinations of non-penetrating snapshots  & & Active-set              \\
        Sparsity induction      & Non-negative FOCUSS     & & Greedy activation and deactivation of dual dofs   \\ 
        Additional assumptions  & Feasible region is convex and parameter-independent & & None                    \\
        \bottomrule
    \end{tabular}
    \caption{Characteristics of the two dictionary-based approximation methods for contact mechanics problems}
	\label{tab:dictionary_methods}
\end{table}

% The inspiration for this approach comes from the reduced modelling of incompressible flow problems, where the velocity field is inherently divergence-free. This property allows the computation of divergence-free snapshots and more importantly, a divergence-free subspace. Consequently, the reduced incompressible flow problem is unconstrained, as all candidates in velocity subspace naturally satisfy the constraints~\cite{Veroy2005,Liberge2010}. Moreover, the orthogonality of the divergence-free velocity field and the irrotational pressure gradient field, results in elimination of pressure term in the reduced balance equation.

The inspiration for this approach comes from the reduced modelling of incompressible flow problems, where the velocity field is inherently divergence-free. Consequently, the reduced incompressible flow problem is unconstrained, as all candidates in velocity subspace naturally satisfy the constraints~\cite{Veroy2005,Liberge2010}. It is not straightforward to extend the same idea to contact mechanics problems, because linear combinations of the penetration-free primal snapshots do not satisfy the non-penetration condition.

    % This approach utilises a monolithic dictionary whose columns are snapshots consisting of both displacement as well as contact pressure. The reasoning for this choice will be evident in the formulation of this approach. Moreover, this method allows the residual to be projected on a low-dimensional space like in the case of unconstrained problems (see \cref{sec:dictionary_resolution}). It also permits the use of sparse regression methods discussed in \cref{sec:sparse_methods}, provided the non-negativity constraint can be enforced. In this work, the non-negative version of the FOCUSS algorithm is used.

    For simplicity, let us first consider a contact problem where the contact pairs are constant. The discrete inequality constraint can be expressed using matrices that do not change with the displacement field. The feasible region $\mathcal{K}$ in the displacement solution space can be defined as:
\begin{align}
    \mathcal{K} := \lbrace \bm{u}  \ | \ \mathbf{C} \bm{u} - \bm{g} \leq 0 \rbrace
\end{align}

\begin{figure}[!htb]
    \centering
    \includegraphics{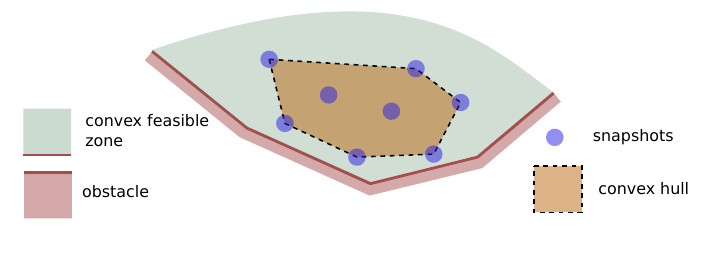}
    \caption{Illustration of non-penetrating property of convex hull. The convex feasible region $\mathcal{K}$ is separated from non-feasible zone by the hyperplanes defined by $\mathbf{C} \bm{u} = \bm{g}$, illustrated schematically by segments. The convex hull of the displacement snapshots is a convex subset of low-dimensional subspace defined by the linear combination of snapshots. It is also a subset of the convex feasible zone $\mathcal{K}$.}
    \label{fig:convex_hull_demo}
\end{figure}

As the operators $\mathbf{C}$ and $\mathbf{g}$ are independent of displacement field $\bm{u}$, it is straightforward to demonstrate the convexity of $\mathcal{K}$. Thus, a convex combination\footnote{A convex combination is a linear combination with non-negative coefficients that sum up to unity} of non-penetrating snapshots would naturally satisfy the inequality constraint. The implication here is that all the candidates in convex hull\footnote{Convex hull of a set of points (or snapshots) is the set containing all convex combinations of the points. In other words, it is the smallest convex set that envelopes the set of points.} of the training set snapshots satisfy the inequality constraint. This notion is visually shown in \cref{fig:convex_hull_demo}, where the convex hull of training set snapshots does not violate the convex non-penetration constraint. This property of convexity can be exploited to build efficient reduced models, if the following hypothesis holds: 

\begin{hypothesis}[Convex Subset hypothesis]\label{hyp:lowdim_convex_subset}
    Given a parametrized high-dimensional inequality-constrained problem with a convex feasible set $\mathcal{K}$, its solutions lie in a low-dimensional convex subset inside $\mathcal{K}$.
\end{hypothesis}

% The Convex Subset hypothesis is an extension of the idea behind the low-rank hypothesis to inequality-constrained problems with convex feasible regions. It allows us to explore a solution set limited to the convex hull defined by the snapshots without worrying about the non-penetration conditions. Note that the convexity of the feasible region assures that any solution in the convex hull of displacement snapshots still lies inside the feasible region, independent of the \cref{hyp:lowdim_convex_subset}. 

To explore the convex hull of the dictionary, the displacement field will be approximated using convex combinations of the primal dictionary $\PrimalDict$:
\begin{align*}
    \bm{u} &\approx \PrimalDict \widehat{\bm{u}} \\
    % \text{s.t. } \sum_{i} \widehat{u}_i &= 1 \\
    \text{s.t. } \bm{1}^T \widehat{\bm{u}} &= 1 \\
    \widehat{\bm{u}} & \geq \bm{0}
\end{align*}

\noindent As the non-penetration constraint will be satisfied naturally by all candidates in the convex hull of dictionary snapshots, a dictionary-based approximation can be computed using coefficients $\widehat{\bm{u}}$ that satisfy equilibrium equations. Therefore, the objective is to minimize the following residual:
\begin{align*}
    \bm{r}(\bm{\alpha}) = \mathbf{K}(\mu) \PrimalDict \bm{\alpha} + \mathbf{C}^T(\mu, \bm{u}) \DualDict \bm{\alpha} - \bm{f}(\mu)
\end{align*}

\noindent Note that the residual has been defined using $\widehat{\bm{u}} = \widehat{\bm{\lambda}} = \bm{\alpha}$, which is equivalent to using a monolithic dictionary, in which each column consists of a displacement snapshot stacked above the corresponding contact pressure snapshot. Usage of a monolithic dictionary not only reduces the number of unknowns but also makes the problem of minimizing the residual $\norm{\bm{r}(\bm{\alpha})}{}$ well-defined. Moreover, the residual can be projected on a low-dimensional subspace of the $\mathtt{colsp}(\PrimalDict)$, described by the matrix $\mathbf{B}$ containing its basis, while also adding sparsity constraints to the unknown $\bm{\alpha}$, like in~\eqref{eq:sparse_res_randomized}. Thus, the problem can be stated as:

\begin{subequations}
    \begin{align}
        \min \ & \norm{\vphantom{\widehat{}} \bm{\alpha}}{p} \\
        \text{s.t. } & \norm{\mathbf{B}^T \bm{r}(\bm{\alpha})}{} < \varepsilon \\
            &  \bm{1}^T\bm{\alpha} = 1 \label{eq:convex_sum} \\
            &  \bm{\alpha} \geq \bm{0}
    \end{align} \label{eq:convex_comb_equilibrium}
\end{subequations}

% As the constraint~\eqref{eq:kkt_nonpenetration} is not imposed directly in this approach, the non-penetration condition is satisfied only because of the convexity of the feasible region. Thus, the equality condition $\mathbf{C}_{\ActiveSet}(\mu,\bm{u}) \bm{u} = \bm{g}_{\ActiveSet}(\mu,\bm{u})$ for active constraints is also not imposed explicitly. Consequently, the complementary slackness KKT condition~\eqref{eq:kkt_comp_slack} may not be satisfied exactly. This will be seen in the upcoming numerical examples. \\

The nnFOCUSS algorithm~\cite{Lauzeral2019,Manucci2022}, based on the FOCUSS algorithm~\cite{Gorodnitsky1997}, is a weighted-norm minimization technique that progressively induces sparsity to the solution(see \cref{sec:sparse_methods_list} for details). The algorithm takes three inputs: the vector $\bm{x}$ that will be approximated, the dictionary $\mathbf{D}$ that serves as the approximation space, and a non-sparse non-negative initial guess $\bm{\alpha}^0$ of approximation coefficients. It can be used to solve the convex hull approximation problem~\eqref{eq:convex_comb_equilibrium}, by plugging the inputs defined as:
\begin{align*}
    \mathbf{D} =  \begin{bmatrix}
        \mathbf{B}^T  \mathbf{K}_{\texttt{mono}}\\
        \bm{1}^T
    \end{bmatrix} , \ \ \ \  % \\
    \bm{x} = \begin{bmatrix}
        \mathbf{B}^T \bm{f}(\mu) \\
        1
    \end{bmatrix} , \ \ \ \  % \\
    \bm{\alpha}^0 = \mathtt{nnls}(\mathbf{D}, \bm{x})
\end{align*}
\noindent where $\mathbf{K}_{\texttt{mono}}$ is the operator for the monolithic residual $\bm{r}(\bm{\alpha})$:
\begin{align*}
    \mathbf{K}_{\texttt{mono}} = \mathbf{K}(\mu) \PrimalDict  + \mathbf{C}^T(\mu, \bm{u}) \DualDict 
\end{align*}
\texttt{nnls} indicates the non-negative least squares. An open-source implementation of \texttt{nnls} provided in the SciPy library~\cite{scipy} is used.

A test is proposed which checks if \cref{hyp:lowdim_convex_subset} hold within a reasonable error using the so-called leave-one-out approach on the training set snapshots. The \cref{hyp:lowdim_convex_subset} can be verified by approximating the left-out snapshot with the convex hull of the rest of the snapshots using least-square criteria. The approximation error, referred to here as the convex hull least square (CHLS) error, can be used to assess the existence of a hypothetical low-dimensional convex set. Each dictionary element $\bm{d}$ is approximated using the convex hull of $\widetilde{\PrimalDict}$ and the error is recorded:
        \begin{gather}
            \begin{aligned}
                \epsilon_{\mathtt{CHLS}}[k] = \min_{\bm{\alpha}} \frac{\norm{\widetilde{\PrimalDict} \bm{\alpha} - \bm{d}}{2}}{\norm{\bm{d}}{2}} \\
                \text{s.t.}  \bm{1}^T \bm{\alpha} = 1
                % \label{eq:test_lowdim_chull}
            \end{aligned}
        \end{gather}
        A high value of $\epsilon_{\mathtt{CHLS}}$ might occur due to a lack of low-dimensional convexity or due to scattered training set data.

\subsection{Illustrative example: elastic rope-obstacle problem}\label{sec:convex_rope}
The method is demonstrated on a problem where the convexity of feasible space is assured. The elastic rope-obstacle problem from references~\cite{Haasdonk2012,Balajewicz2016,Bader2016} is considered, where a 1D elastic rope, fixed at the two ends comes into contact with a fixed obstacle (\cref{fig:membrane_snaps_g1}). The elasticity $\nu(x, \gamma)$ of the rope is a function of parameter $\gamma$.
\begin{gather}
	\begin{aligned}
        &\nu(x, \gamma) \nabla^2 u(x) = f \ \ \text{ on } x \in [0,1] \\
		&u(0) = u(1) = 0 \\
		&u(x) \geq -0.2 (\sin(\pi x) - \sin(3 \pi x)) -0.5 \\
        \text{where} &  \\
        &\nu(x, \gamma) = \left \{  
            \begin{matrix} \gamma  & \text{if } x < 0.5 \\
                                30      & \text{otherwise}
            \end{matrix} \right . \ \ \text{ on } \gamma \in [10,50]
	\end{aligned}
	\label{eq:membrane_elastic_parameter}
\end{gather}

\noindent $
\begin{aligned}
    &\text{Training set:} && \gamma \in \mathcal{P}_{\texttt{tr}}\hspace{0.5em} = \{10, 15, 20 \dots 45, 50\}\\
    &\text{Validation set:} && \gamma \in \mathcal{P}_{\texttt{val}} = \{12.5, 17.5, 22.5 \dots  42.5, 47.5\}
\end{aligned}
$

\begin{figure}[htpb!]
    \centering
    \includegraphics[width=0.5\linewidth]{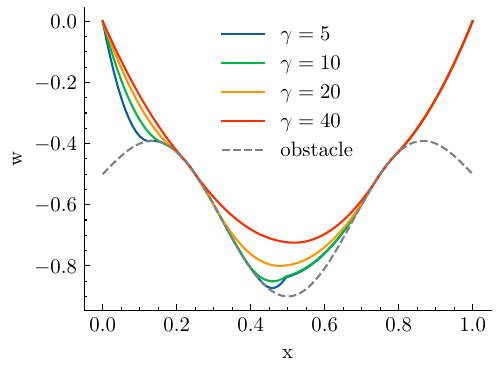}
    \caption{Elastic rope-obstacle problem: Snapshots of the elastic rope fixed at the two ends, while it deforms under a uniform load and establishes contact with the obstacle. Snapshots for various values of $\gamma$ are shown.}
    \label{fig:membrane_snaps_g1}
\end{figure}

The monolithic dictionary with 9 snapshots corresponding to the training set defined above is computed in the offline stage. The test results for \cref{hyp:lowdim_convex_subset} are shown in \cref{fig:membrane_lowdim_convex}, indicating that the hypothesis holds to a reasonable level.

\begin{figure}[htpb!]
    \centering
    \includegraphics[width=0.5\linewidth]{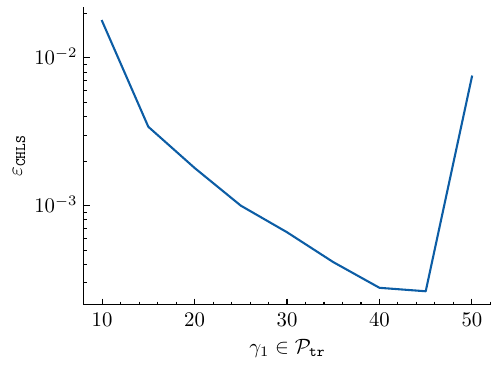}
    \caption{Testing the \cref{hyp:lowdim_convex_subset} for the given training set of elastic rope-obstacle problem~\eqref{eq:membrane_elastic_parameter}}
    \label{fig:membrane_lowdim_convex}
\end{figure}

The low-rank matrix $\mathbf{B}$ on which the residual will be projected is computed using left singular vectors of the truncated SVD $\mathbf{B} \gets \text{svd}(\PrimalDict, \delta)$, with truncation tolerance $\delta = 10^{-7}$. The reconstruction errors and resultant sparsity of the solution in the training set and a validation set are shown in \cref{fig:convex_membrane_gamma1}. The training set is reconstructed within near numerical precision. In the validation set, the algorithm chooses the training set points that are nearest to $\gamma$ to create the best possible reconstruction with the given dictionary. 

As discussed earlier, the non-penetration condition and complementary slack are not explicitly imposed. However, as the feasible region is convex, one would expect that at least the non-penetration condition would be satisfied exactly. Nonetheless, as seen in the \cref{tab:convex_membrane_kkt_values}, reconstructions in validation set show a finite but small penetration and complementary slack. These violations appear because the condition~\eqref{eq:convex_sum} is satisfied only to a certain precision. The violation to convex condition, evaluated as $|\bm{1}^T \bm{\alpha}-1|$ in \cref{tab:convex_membrane_kkt_values}, is of the similar order as penetration and complementary slack.

\begin{table}[!htb]
	\setlength\extrarowheight{2pt}
    \small
	\centering
    \begin{tabular}{lcc}
        \toprule
        &  $\gamma \in \mathcal{P}_{\texttt{tr}}$   & $\gamma \in \mathcal{P}_{\texttt{val}}$  \\
        \midrule
        $|\bm{1}^T \bm{\alpha}-1|$ & $\mathcal{O}(10^{-13})$ & $\mathcal{O}(10^{-3})$ \\
        Penetration & $\mathcal{O}(10^{-16})$ & $\mathcal{O}(10^{-4})$ \\
        Complementary Slack & $\mathcal{O}(10^{-15}) $ & $\mathcal{O}(10^{-3})$\\
        \bottomrule
    \end{tabular}
    \caption{Observed values of various compliances for the elastic rope-obstacle problem with parameter $\gamma$.}
	\label{tab:convex_membrane_kkt_values}
\end{table}

\begin{figure}[!htb]
    \centering
    \begin{subfigure}[t]{0.45\linewidth}
        \includegraphics[width=1\linewidth]{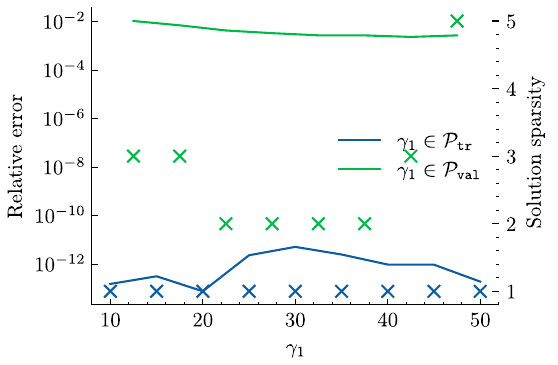}
        \caption{Displacement}
    \end{subfigure}
    \quad
	\begin{subfigure}[t]{0.45\linewidth}
        \includegraphics[width=1\linewidth]{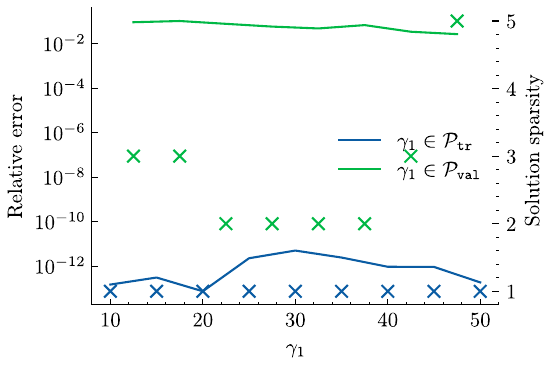}
        \caption{Contact Pressure}
    \end{subfigure}
    \caption{Convex hull reconstruction errors for elastic rope-obstacle problem with parameter $\gamma$ using $\mathbf{B} \gets \text{svd}(\PrimalDict, \delta=10^{-7})$. Crosses $\times$ indicate achieved sparsity (w.r.t.\ right y-axis). Primal and dual sparsity are naturally equal as the dictionary is monolithic.}
    \label{fig:convex_membrane_gamma1}
\end{figure}

\section{Spurious effects of POD trucation}
\label{sec:tau_delta}
In the greedy active-set method, the snapshots of the dual dictionary are selected based on the violations of the non-penetration constraints. The dictionary vector with the highest correlation with the current state of penetration is chosen. However, instead of zero penetration, a small value of penetration $\tau$ is allowed (see~\eqref{eq:sparse_greedy_kkt_nonpen} and~\eqref{eq:sparse_enrich_concept}). This is done to avoid spurious selection of dictionary elements, which happens due to a truncated primal basis. Such a reconstruction example of Hertz problem for parameter value $d=0.25$ is shown in \cref{fig:hertz_reconstruction_spurious}, where a spurious peak in the contact pressure is evident. In this case, a dual dictionary of size 30 and primal basis truncated at $10^{-6}$ are used. The algorithm selected the appropriate snapshots from the dual dictionary, but also selected some unnecessary snapshots: the first ($\DualDict[1]$) and twelfth ($\DualDict[12]$) dictionary vectors. Importantly, the $\DualDict[1]$ vector contributes to large error in reconstruction as evident in the figure. The false selection of this snapshot is linked to the truncation of primal basis, and therefore, can be avoided by setting the value of $\tau$ same as $\delta$.

\begin{figure}[!htb]
	\centering
    \includegraphics[width=0.5\linewidth]{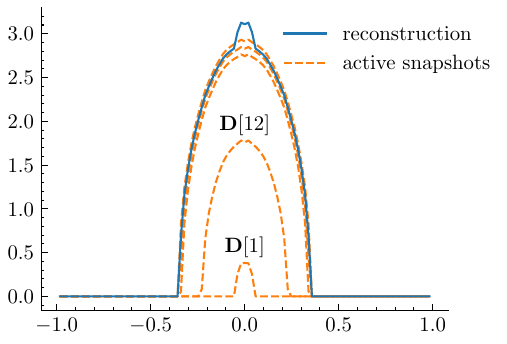}
    \caption
    {Reconstruction instance of Hertz problem with $d=0.25$ using $\# \mathcal{P}_{\texttt{tr}} = 30$ and primal basis truncated at $\delta=10^{-6}$, where few dictionary snapshots were spuriously selected by the greedy active-set algorithm if non-penetration condition is applied ``hardly''. The falsely selected snapshots are indicated by their index in the dual dictionary. The snapshot which contributes to most of the error is $\DualDict[1]$, corresponds to loading parameter $d=0.01$}
	\label{fig:hertz_reconstruction_spurious}
\end{figure}
%%fakesection about the authors
\ifx\fmt\undefined
\section*{About the authors}
\subsection*{Kiran Sagar Kollepara, Jos\'e V. Aguado, Luisa Silva, Domenico Borzacchiello}
  Institut de Recherche en Génie Civil et Mécanique (GeM) at {E}cole Centrale de Nantes\\ 1 rue de la No{\"e}, BP 92101, 44321 Nantes cedex 3, France \\
  e-mail: \{kiran-sagar.kollepara, jose.aguado-lopez, luisa.rocha-da-silva, domenico.borzacchiello\}@ec-nantes.fr\\ 
   web: http://gem.ec-nantes.fr

\subsection*{Yves Le Guennec}
  IRT Jules Verne\\ 1 Mail des 20 000 Lieues, 44340 Bouguenais, France \\ 
   web: http://irt-jules-verne.fr

\fi
\end{document}